# A DYNAMIC LOOK-AHEAD MONTE CARLO ALGORITHM FOR PRICING BERMUDAN OPTIONS

By Daniel Egloff, Michael Kohler and Nebojsa Todorovic

*Zurich Cantonal Bank, University of Saarbrücken and
University of Saarbrücken*

Under the assumption of no-arbitrage, the pricing of American and Bermudan options can be casted into optimal stopping problems. We propose a new adaptive simulation based algorithm for the numerical solution of optimal stopping problems in discrete time. Our approach is to recursively compute the so-called continuation values. They are defined as regression functions of the cash flow, which would occur over a series of subsequent time periods, if the approximated optimal exercise strategy is applied. We use nonparametric least squares regression estimates to approximate the continuation values from a set of sample paths which we simulate from the underlying stochastic process. The parameters of the regression estimates and the regression problems are chosen in a data-dependent manner. We present results concerning the consistency and rate of convergence of the new algorithm. Finally, we illustrate its performance by pricing high-dimensional Bermudan basket options with strangle-spread payoff based on the average of the underlying assets.

**1. Introduction.** Many financial contracts allow for early exercise before expiry. Most of the exchange traded option contracts are of the American type which allows the holder to choose any exercise date before expiry, or the Bermudan with exercise dates restricted to a predefined discrete set of dates. Mortgages have embedded prepayment options such that the mortgage can be amortized or repaid. Also, life insurance contracts may allow for early surrender. In this paper we are interested in pricing options with early exercise features. It is well known that in complete and arbitrage free markets the price of a derivative security can be represented as an expected value with respect to the so-called martingale measure; see, for instance, [16]. Furthermore, the price of an American option with maturity $T$ is given by the value of the optimal stopping problem

$$(1.1) \qquad V_0 = \sup_{\tau \in \mathcal{T}_{[0,T]}} \mathbf{E}\{d_{0,\tau} f_\tau(X_\tau)\},$$

where $f_t$ is a nonnegative payoff function, $X_t$ is a stochastic process, which models the relevant risk factors, $\mathcal{T}_{[0,T]}$ is the class of all stopping times with values in $[0, T]$, and $d_{s,t}$ are nonnegative $\mathcal{F}((X_u)_{s \le u \le t})$-measurable discount factors satisfying $d_{0,t} = d_{0,s} \cdot d_{s,t}$ for $s < t$. In practice, the process $X_t$ is often a geometric Brownian motion, as, for instance, in the celebrated Black–Scholes setting. A more general class of models is obtained with diffusions, jump-diffusion processes or nonparametric time series models. The model parameters are usually calibrated to observed time series data.

The first step in addressing the numerical solution of (1.1) is to pass from continuous time to discrete time, which means in financial terms to approximate the American option by a Bermudan option. The convergence of the discrete time approximations to the continuous time optimal stopping problem is considered in [18] for the Markovian case but also in the abstract setting of general stochastic processes.

For simplicity, we restrict ourselves directly to a discrete time scale and consider exclusively Bermudan options. In analogy to (1.1), the price of a Bermudan option is the value of the discrete time optimal stopping problem

$$(1.2) \qquad V_0 = \sup_{\tau \in \mathcal{T}(0,\ldots,T)} \mathbf{E}\{d_{0,\tau} f_\tau(X_\tau)\},$$

where $X_0, X_1, \ldots, X_T$ is now a discrete time stochastic process, and $\mathcal{T}(0, \ldots, T)$ is the class of all $\{0, \ldots, T\}$-stopping times. For additional theoretical background on valuating Bermudan options, we refer to [25].

In the sequel we assume that $X_0, X_1, \ldots, X_T$ is a $[-A, A]^d$-valued Markov process recording all necessary information about financial variables including prices of the underlying assets as well as additional risk factors driving stochastic volatility or stochastic interest rates. We also assume that the law of $X_0, \ldots, X_T$ is known such that we can draw random sample paths as well as partial sample paths $X_t, \ldots, X_T$ for arbitrary starting values of $X_t$. Neither the Markov property nor the form of the payoff as a function of the state of $X_t$ is restrictive and can always be achieved by including supplementary variables. For instance, in the case of an Asian option we add the running mean as an additional variable into $X_t$. Because the diffusion, jump-diffusion or time series models, which appear in practical applications, lead to unbounded stochastic processes for the underlying state variables $X_t$, they must be suitably localized to a bounded set $[-A, A]^d$.

The boundedness assumption $X_t \in [-A, A]^d$ then allows us to estimate the price of the Bermudan option from samples of polynomial size in the number of free parameters. This is in contrast to Glasserman and Yu [13]. Their work does not impose a boundedness assumption on the underlying process and shows that for arithmetic and geometric Brownian motions, the sample size must grow exponentially in the number of free parameters in order to retain a convergent estimator.

The computation of (1.2) requires the determination of an optimal stopping rule $\tau^* \in \mathcal{T}(0, \ldots, T)$ which satisfies

$$(1.3) \qquad V_0 = \mathbf{E}\{d_{0,\tau^*} f_{\tau^*}(X_{\tau^*})\}.$$

Let

$$(1.4) \qquad q_t(x) = \sup_{\tau \in \mathcal{T}(t+1, \ldots, T)} \mathbf{E}\{d_{t,\tau} f_\tau(X_\tau)|X_t = x\}$$

be the so-called continuation value describing the value of the option at time $t$ given $X_t = x$ and subject to the constraint of holding the option at time $t$ rather than exercising it. The general theory of optimal stopping for Markov processes (see, e.g., [5, 11, 22, 26]) implies that

$$\tau^* = \inf\{s \geq 1 : q_s(X_s) \leq f_s(X_s)\}$$

is an optimal stopping time, that is, $\tau^*$ satisfies (1.3). Therefore, computing the continuation values (1.4) solves the optimal stopping problem (1.2).

Explicit solutions of (1.2) do not exist, except in very rare cases, but there are a variety of numerical procedures to solve optimal stopping problems, each with its strength and weaknesses. In this paper we study a concrete simulation algorithm. The first attempts to use simulation are [2, 3, 28]. Longstaff and Schwartz [21] introduce a new algorithm for Bermudan options in discrete time. It combines Monte Carlo simulation with multivariate function approximation. Tsitsiklis and Van Roy [29] independently propose an alternative parametric approximation algorithm using stochastic approximation to derive the weights of the approximation. Both algorithms approximate the value function or the early exercise rule and therefore provide a lower bound for the true optimal stopping value. Upper bounds based on the dual problem are derived in [15, 23]. More details and further references can be found in [4] and [12]. The article [19] compares several Monte Carlo approaches empirically.

In this paper we enhance the approach of [21] and its generalization presented in [10]. We construct estimates $\hat{q}_t$ of $q_t$ and approximate the optimal stopping rule $\tau^*$ by

$$(1.5) \qquad \hat{\tau} = \inf\{s \geq 1 : \hat{q}_s(X_s) \leq f_s(X_s)\}.$$

Then, a Monte Carlo estimate of

$$(1.6) \qquad \mathbf{E}\{d_{0,\hat{\tau}} f_{\hat{\tau}}(X_{\hat{\tau}})\}$$

provides a lower bound for the price $V_0$ of the Bermudan option.

To this end, we represent $q_t$ as a regression function of a distribution $(X_t, Y_t)$, where $Y_t$ depends on the partial sample path $X_{t+1}, \ldots, X_{t+w+1}$ and $q_{t+1}, \ldots, q_{t+w+1}$ for some tunable parameter $w \in \{0, 1, \ldots, T - t - 1\}$. This distribution will in turn be approximated by $(X_t, \hat{Y}_t)$, where $\hat{Y}_t$ depends on $X_{t+1}, \ldots, X_{t+w+1}$ and $\hat{q}_{t+1}, \ldots, \hat{q}_{t+w+1}$. We construct an estimate $\hat{q}_t$ of $q_t$ with nonparametric regression techniques applied to a Monte Carlo sample of the distribution $(X_t, \hat{Y}_t)$ and use this estimate together with $\hat{q}_{t+1}, \ldots, \hat{q}_{t+w}$ to compute recursive estimates of $q_{t-1}, \ldots, q_0$. Our algorithm is adaptive in the sense that all parameters of the estimates and the parameter $w$ of the distribution of $(X_t, Y_t)$ are chosen in a data dependent manner.

We proceed as follows. In Section 2 we describe in detail the connection between discrete time optimal stopping problems and recursive regression. The dynamic look-ahead Monte Carlo algorithm for solving optimal stopping problems is introduced in Section 3. The main theoretical results, including the consistency and the rate of convergence of the algorithm, are presented in Section 4. The finite sample properties of the proposed algorithm are illustrated in Section 5 with a simulation study. Section 6 contains the proofs.

## 2. Discrete time optimal stopping and recursive regression.

Let $\mathbf{X} = (X_t)_{t=0,\ldots,T}$ be a discrete time Markov process with values in $\mathbb{R}^d$, $\mu_t$ the law induced by $X_t$ on $\mathbb{R}^d$, and $\mathbb{F} = (\mathcal{F}_t)$ be the induced filtration where

$$(2.1) \qquad \mathcal{F}_t = \mathcal{F}(X_0, \ldots, X_t) = \bigvee_{s \le t} \sigma(X_s)$$

is the sigma algebra generated by the random variables $\{X_s | s \le t\}$. The solution of the discrete time optimal stopping problem for nonnegative reward or payoff functions $f_t$ is given by the value function

$$(2.2) \qquad v_t(x) = \sup_{\tau \in \mathcal{T}(t,\ldots,T)} \mathbf{E}[f_\tau(X_\tau) | X_t = x].$$

The supremum runs over the class $\mathcal{T}(t, \ldots, T)$ of all $\mathbb{F}$-stopping times with values in $\{t, \ldots, T\}$. By definition, each $\tau \in \mathcal{T}(t, \ldots, T)$ satisfies $\{\tau = k\} \in \mathcal{F}(X_0, \ldots, X_k)$ for $k \in \{t, \ldots, T\}$. Here and in the sequel we assume for notational simplicity that $f_t$ contains already the discount factor occurring in (1.2). Once the value function has been determined, the smallest optimal stopping time as of time $t$ can be derived as

$$(2.3) \qquad \tau_t^* = \inf\{s \ge t | v_s(X_s) \le f_s(X_s)\}.$$

The optimal stopping problem can also be characterized in terms of the so-called continuation value, which is given by

$$(2.4) \qquad q_t(x) = \sup_{\tau \in \mathcal{T}(t+1,\ldots,T)} \mathbf{E}[f_\tau(X_\tau) | X_t = x] = \mathbf{E}[f_{\tau_{t+1}^*}(X_{\tau_{t+1}^*}) | X_t = x]$$

for $t \leq T - 1$ and set to $q_T = 0$ at maturity $T$. The value function and the continuation value are related by

$$(2.5) \quad v_t(X_t) = \max(f_t(X_t), q_t(X_t)), \qquad q_t(X_t) = \mathbf{E}[v_{t+1}(X_{t+1})|X_t].$$

From now on we primarily consider $q_t$. The continuation value satisfies the dynamic programming equations

$$(2.6) \quad \begin{aligned} &q_T(x) = 0, \\ &q_t(x) = \mathbf{E}[\max(f_{t+1}(X_{t+1}), q_{t+1}(X_{t+1}))|X_t = x]. \end{aligned}$$

The recursion for the optimal stopping rules is given by

$$(2.7) \quad \begin{aligned} &\tau_T^* = T, \\ &\tau_t^* = t \mathbf{1}_{\{q_t(X_t) \leq f_t(X_t)\}} + \tau_{t+1}^* \mathbf{1}_{\{q_t(X_t) > f_t(X_t)\}}. \end{aligned}$$

The dynamic programming equations (2.6) show that the optimal stopping problem in discrete time is essentially equivalent to a series of regression problems. Equation (2.4) provides a different regression representation of the continuation value, once the optimal stopping rule of the next future period is known. These representations are extreme cases, as we will explain in the following. For $h_t \in L_1(\mu_t)$ with $h_T = f_T$, we define on $\mathbb{R}^{(w+1)d} = \times_{w+1} \mathbb{R}^d$ the function

$$(2.8) \quad \begin{aligned} &\vartheta_{t:w}(f, h_t, \ldots, h_{t+w})(x_t, \ldots, x_{t+w}) \\ &= \sum_{s=t}^{t+w} f_s(x_s) \mathbf{1}_{\{f_s(x_s) - h_s(x_s) \geq 0\}} \prod_{r=t}^{s-1} \mathbf{1}_{\{f_r(x_r) - h_r(x_r) < 0\}} \\ &\quad + h_{t+w}(x_{t+w}) \prod_{r=t}^{t+w} \mathbf{1}_{\{f_r(x_r) - h_r(x_r) < 0\}}, \end{aligned}$$

where we follow the convention that the product over an empty index set is equal to one. In the following, to reduce notational overhead, we simply write

$$(2.9) \quad \vartheta_{t:w}(f, h) = \vartheta_{t:w}(f, h_t, \ldots, h_{t+w}),$$

thereby implicitly assuming that $\vartheta_{t:w}(f, h)$ is solely depending on $h_t, \ldots, h_{t+w}$.

In a financial context the function $\vartheta_{t:w}(f, h)$ has a natural interpretation as the future payoff we would get by holding the Bermudan option for at most $w$ periods, applying the stopping rule $\tau_t(h) \wedge (t + w)$ which is defined recursively by

$$(2.10) \quad \begin{aligned} &\tau_T(h) = T, \\ &\tau_t(h) = t \mathbf{1}_{\{f_t(X_t) - h_t(X_t) \geq 0\}} + \tau_{t+1}(h) \mathbf{1}_{\{f_t(X_t) - h_t(X_t) < 0\}}, \end{aligned}$$

and selling the option at time $t + w$ for the price $h_{t+w}(X_{t+w})$, if it is not exercised before.

We now come back to the generalization of the regression representations (2.4) and (2.6). First note that $\max(f_{t+1}, q_{t+1}) = \vartheta_{t+1:0}(f, q)$ and, therefore,

$$(2.11) \qquad q_t(x) = \mathbf{E}[\vartheta_{t+1:0}(f, q)(X_{t+1})|X_t = x].$$

On the other hand, the recursive formula (2.7) for the optimal stopping rule $\tau_t^*$ shows that

$$f_{\tau_{t+1}^*}(X_{\tau_{t+1}^*}) = f_{\tau_{t+1}(q)}(X_{\tau_{t+1}(q)}) = \vartheta_{t+1:T-t-1}(f, q)(X_{t+1}, \dots, X_T),$$

such that we also have [cf. (2.4)]

$$(2.12) \qquad q_t(x) = \mathbf{E}[\vartheta_{t+1:T-t-1}(f, q)(X_{t+1}, \dots, X_T)|X_t = x].$$

More generally, we have for any $0 \le w \le T - t - 1$ the representation

$$(2.13) \qquad q_t(x) = \mathbf{E}[\vartheta_{t+1:w}(f, q)(X_{t+1}, \dots, X_{t+w+1})|X_t = x].$$

To prove (2.13), we start with

$$
\begin{aligned}
(2.14) \quad q_t(X_t) &= \mathbf{E}[\max(f_{t+1}(X_{t+1}), q_{t+1}(X_{t+1}))|X_t] \\
&= \mathbf{E}[f_{t+1}(X_{t+1})1_{\{f_{t+1}(X_{t+1}) - h_{t+1}(X_{t+1}) \ge 0\}} \\
&\quad + q_{t+1}(X_{t+1})1_{\{f_{t+1}(X_{t+1}) - h_{t+1}(X_{t+1}) < 0\}}|\mathcal{F}_t],
\end{aligned}
$$

where we have used the Markov property in the second equality. Then we expand $q_{t+1}(X_{t+1})$ in (2.14) by

$$
\begin{aligned}
\mathbf{E}[&f_{t+2}(X_{t+2})1_{\{f_{t+2}(X_{t+2}) - h_{t+2}(X_{t+2}) \ge 0\}} \\
&+ q_{t+2}(X_{t+2})1_{\{f_{t+2}(X_{t+2}) - h_{t+2}(X_{t+2}) < 0\}}|\mathcal{F}_{t+1}]
\end{aligned}
$$

and proceed recursively up to $t + w + 1$. Equation (2.13) follows from the projection property $\mathbf{E}[\mathbf{E}[\cdot|\mathcal{F}_{t+1}]|\mathcal{F}_t] = \mathbf{E}[\cdot|\mathcal{F}_t]$ of conditional expectations and by another application of the Markov property.

**3. Monte Carlo algorithms for optimal stopping.** Equation (2.13) shows that the continuation value $q_t$ at time $t$ can be obtained as the regression function of $\vartheta_{t+1:w}(f, q)$ for some $0 \le w \le T - t - 1$. Least squares Monte Carlo methods pioneered by [21], and extended in [10] to arbitrary $w$, recursively estimate the regression function $q_t$ from independent sample paths of the underlying Markov process $X_t$. Let

$$(3.1) \qquad X_{t+1:w} = (X_{t+1}, \dots, X_{t+w+1})$$

be the partial sample path of length $w$ starting at $t + 1$. When it comes to estimation of the continuation value $q_t$, these algorithms use the previously determined estimates $\hat{q}_{t+1}, \dots, \hat{q}_{t+w+1}$ for $q_{t+1}, \dots, q_{t+w+1}$ to construct

$$(3.2) \quad \hat{Y}_t = \vartheta_{t+1:w}(f, \hat{q})(X_{t+1:w}) = \vartheta_{t+1:w}(f, \hat{q}_{t+1}, \dots, \hat{q}_{t+w+1})(X_{t+1:w}),$$

which takes the role of the dependent variable of the regression problem for time step $t$. The random variable $\hat{Y}_t$ is an estimate of the unknown optimal reward

$$(3.3) \qquad Y_t = \vartheta_{t+1:w}(f,q)(X_{t+1:w}) = \vartheta_{t+1:w}(f,q_{t+1},\ldots,q_{t+w+1})(X_{t+1:w}).$$

Given independent sample paths

$$(3.4) \qquad \mathbf{X}_i = (X_{i,t})_{t=0,\ldots,T}, \qquad i = 1,\ldots,n,$$

of the underlying Markov process $\mathbf{X}$, the least squares estimate of $q_t$ is obtained as

$$(3.5) \qquad \hat{q}_{n,t} = \arg\min_{h \in \mathcal{H}_{n,t}} \frac{1}{n} \sum_{i=1}^{n} |h(X_{i,t}) - \hat{Y}_{i,t}|^2,$$

where

$$(3.6) \qquad \hat{Y}_{i,t} = \vartheta_{t+1:w}(f,\hat{q})(X_{i,t+1:w}), \qquad X_{i,t+1:w} = (X_{i,t+1},\ldots,X_{i,t+w+1})$$

and $\mathcal{H}_{n,t}$ is a set of functions $h : \mathbb{R}^d \to \mathbb{R}$.

With $w = 0$, the above algorithm corresponds to the Tsitsiklis–Van Roy algorithm [29], while $w = T - t - 1$ has been proposed in [21]. The idea of using an intermediate value $w \in \{0,1,\ldots,T-t-1\}$ in order to "interpolate" between these two algorithms has been introduced in [10]. A further contribution of [10] is the consistency and the rate of convergence of the above algorithm for fixed $w$ and fixed convex and uniformly bounded function spaces $\mathcal{H}_{n,t}$, without imposing any distributional assumptions on the underlying process $X_t$.

The boundedness assumption on $\mathcal{H}_{n,t}$ makes the computation of the least squares estimate in (3.5) difficult because it leads to constrained optimization problems; see, for instance, [14], Section 10.1. In addition, the convexity assumption excludes promising choices like spaces of polynomial splines with free knots or spaces of artificial neural networks, which require restrictions on the number of knots or the number of hidden neurons, respectively, to control the "complexity" of the function spaces. The resulting function spaces violate the convexity assumptions. Taking the convex hull instead is not an option because it would lead to function classes with a complexity that is much too high. Furthermore, in view of applications, it is desirable to choose parameters of the functions spaces and also the parameter $w$ of the underlying regression problems data dependent. In this paper we modify the above algorithm such that this is possible. For simplicity, we restrict ourselves to function spaces, which are linear vector spaces, however, it is straightforward to derive similar results for spaces of polynomial splines with free knots or spaces of artificial neural networks.

The main problem in analyzing the estimates $\hat{q}_{n,t}$ is the control of the error propagation, that is, to answer the question how the errors of $\hat{q}_{n,t+1}, \ldots,$ $\hat{q}_{n,t+w+1}$ influence the error of $\hat{q}_{n,t}$. At this stage Egloff [10] uses the convexity of $\mathcal{H}_{n,t}$ to bound the $L_2$-error in terms of the approximation error and a sample error derived from a suitably centered loss function. The difficulty for obtaining error estimates comes from the fact that $\hat{q}_{t+1},\ldots,\hat{q}_{t+w+1}$ depend on a single set of sample paths (3.4) and are thus dependent. Clément, Lamberton and Protter [6] face the same difficulty while deriving a central limit theorem for the Longstaff–Schwartz algorithm with linear approximation.

In the sequel we use a trick to simplify the analysis of the error propagation. Instead of using the partial sample path $X_{t+1:w}$ of our training data again, which we used in part already in the construction of the estimates $\hat{q}_{n,t+1},\ldots,\hat{q}_{n,t+w+1}$, we generate new data $X_{t+1:w}^{t,\text{new}}$ for $X_{t+1:w}$ which are conditionally independent from all previously used data of time $s > t$ given $X_t$ at time point $t$. We then construct samples of the distribution of $(X_t, \hat{Y}_t^{w,\text{new}})$, where

$$\hat{Y}_t^{w,\text{new}} = \vartheta_{t+1:w}(f, \hat{q}_{n,t+1}, \ldots, \hat{q}_{n,t+w+1})(X_{t+1:w}^{t,\text{new}}).$$

Since for $X_t$ given, the random variable $X_{t+1:w}^{t,\text{new}}$ is independent of all previously used data for all time points $s > t$, it is, in particular, independent of the data used in the construction of $\hat{q}_{n,t+1},\ldots,\hat{q}_{n,t+w+1}$. Set

$$q_t^{w,\text{new}}(x) = \mathbf{E}^*\{\hat{Y}_t^{w,\text{new}} | X_t = x\},$$

where in $\mathbf{E}^*\{\cdot | X_t = x\}$ we take the conditional expectation with respect to fixed $X_t = x$ and with all the data fixed which were used in the construction of $\hat{q}_{n,t+1},\ldots,\hat{q}_{n,t+w+1}$. Proposition 6.4 in [10] implies

$$\left\{ \int |q_t^{w,\text{new}}(x) - q_t(x)|^2 \mu_t(dx) \right\}^{1/2}$$

(3.7)

$$\leq \sum_{s=t+1}^{t+w+1} \left\{ \int |\hat{q}_{n,s}(x) - q_s(x)|^2 \mu_s(dx) \right\}^{1/2}.$$

This allows us to control the error propagation. By induction, assume that we have

$$\mathbf{P}\left\{ \int |\hat{q}_{n,s}(x) - q_s(x)|^2 \mu_s(dx) \right.$$

(3.8)

$$> \sum_{r=s}^{T-1} c \cdot \left( \delta_{n,r} + \min_{h \in \mathcal{H}_{n,r}} \int |h(x) - q_r(x)|^2 \mu_r(dx) \right) \right\}$$

$$\to 0 \qquad (n \to \infty)$$

for $s \in \{t+1,\ldots,t+w+1\}$. Assume, in addition, that we are able to show

$$\mathbf{P}\left\{ \int |\hat{q}_{n,t}(x) - q_t^{w,\text{new}}(x)|^2 \mu_t(dx) \right.$$

$$(3.9) \qquad > c \cdot \left( \delta_{n,t} + \min_{h \in \mathcal{H}_{n,t}} \int |h(x) - q_t^{w,\text{new}}(x)|^2 \mu_t(dx) \right) \Bigg\}$$
$$\to 0 \qquad (n \to \infty),$$

which is for suitable $\delta_{n,t}$ (depending on the "complexity" of the function spaces $\mathcal{H}_{n,t}$) a standard rate of convergence result for least squares estimates from a sample of size $n$, where in the sample the response variables are independent given the predictor variables and where the predictor variables are independent; see [30] or [17].

It can be shown that (3.7)–(3.9) imply

$$\mathbf{P} \Bigg\{ \int |\hat{q}_{n,t}(x) - q_t(x)|^2 \mu_t(dx)$$
$$> \bar{c} \cdot \sum_{s=t}^{T-1} \left( \delta_{n,s} + \min_{h \in \mathcal{H}_{n,s}} \int |h(x) - q_s(x)|^2 \mu_s(dx) \right) \Bigg\}$$
$$\to 0 \qquad (n \to \infty).$$

Details concerning related arguments can be found in the proofs of Theorems 4.1 and 4.4 below.

The main difference between our work here and the algorithms used in [21] and [10] is that we generate new data to construct samples of $\hat{Y}_t^{w,\text{new}}$. Therefore, the data used for estimation of $q_t^{w,\text{new}}$ is conditionally independent given the sample of $X_t$, which enables us to conclude (3.9) from standard rate-of-convergence results in nonparametric regression. The generation of the new, independent data is similar to the data generation in the random tree method (see, e.g., Section 8.3 in [12]). However, in contrast to the random tree method, we use nonparametric regression techniques to estimate the regression function, while in the random tree method simple averages are used to estimate the regression function point by point. As a consequence, the number of data points for the random tree method grows exponentially in $T$, while for our method it grows only linearly in $T$.

In the sequel we explain the definition of the estimates in detail. Let $n$ be the number of samples which we generate for our regression estimates, and let $w_{\max} \in \{0, 1, \ldots, T-1\}$ be the maximal look-ahead which we use. We start with generating $n$ independent sample paths

$$\mathbf{X}_i = (X_{i,t})_{t=0,\ldots,T} \qquad (i = 1, \ldots, n)$$

of the underlying Markov process $\mathbf{X}$. Then we set

$$\hat{q}_T = \hat{q}_{n,T} = 0$$

and construct successively estimates of $q_{T-1}, \ldots, q_0$ as follows: Fix $t \in \{0, 1, \ldots, T-1\}$ and assume that estimates $\hat{q}_{n,t+1}, \ldots, \hat{q}_{n,T-1}$ of $q_{t+1}, \ldots, q_{T-1}$ are already constructed. Let

$$w_{\max}(t) = \min\{w_{\max}, T - t - 1\}$$

be the maximal look-ahead of time period $t$. Generate independent sample paths

$$X_{i,t:w_{\max}(t)+1}^{t,\text{new}} = (X_{i,s}^{t,\text{new}})_{s=t,\ldots,t+w_{\max}(t)+1} \qquad (i = 1, \ldots, n)$$

starting at $X_{i,t}^{t,\text{new}} = X_{i,t}$ for every $i \in \{1, \ldots, n\}$ such that, for all $i$, the partial sample paths

$$(3.10) \qquad\qquad\qquad X_{i,t:w_{\max}(t)+1}^{t,\text{new}}$$

have the same distribution as $X_{i,t:w_{\max}(t)+1}$, and such that, given $X_{1,t}, \ldots, X_{n,t}$, this data is independent of all previously generated data points for all time points $s > t$. Define

$$\hat{Y}_{i,t}^{w,\text{new}} = \vartheta_{t+1:w}(f, \hat{q}_{n,t+1}, \ldots, \hat{q}_{n,t+w+1})(X_{i,t+1}^{t,\text{new}}, \ldots, X_{i,t+w+1}^{t,\text{new}})$$

for every $w \in \{0, \ldots, w_{\max}(t)\}$ and apply a nonparametric least squares estimate to the data

$$(3.11) \qquad\qquad\qquad ((X_{i,t}, \hat{Y}_{i,t}^{w,\text{new}}))_{i=1,\ldots,n}$$

to construct estimates $\hat{q}_{n,t}^w$ of $q_t$. The final step is to choose

$$\hat{w}_t \in \{0, 1, \ldots, w_{\max}(t)\}.$$

The resulting estimator for $q_t$ is then given by

$$(3.12) \qquad\qquad\qquad \hat{q}_{n,t} = \hat{q}_{n,t}^{\hat{w}_t}.$$

Next, we explain in detail how to define the nonparametric least squares estimates applied to the data (3.11) and how to select $\hat{w}_t$ in a data dependent way. To this end, we split the sample in three parts: a learning sample of size $n_l$, a testing sample of size $n_t$ and a validation sample of size $n_v$, where $n = n_l + n_t + n_v$. Furthermore, we assume that we are given a finite set $\mathcal{P}_n$ of parameters and for each $p \in \mathcal{P}_n$, a set $\mathcal{H}_{n,p}$ of functions $h : \mathbb{R}^d \to \mathbb{R}$.

For $w$ fixed, we first define $\hat{q}_{n,t}^w$. For every $p \in \mathcal{P}_n$, let

$$(3.13) \qquad \tilde{q}_{n,t}^{w,p}(\cdot) = \underset{h \in \mathcal{H}_{n,p}}{\arg\min} \frac{1}{n_l} \sum_{i=1}^{n_l} |h(X_{i,t}) - \hat{Y}_{i,t}^{w,\text{new}}|^2$$

be the least squares estimate of $q_t^{w,\text{new}}$ in $\mathcal{H}_{n,p}$, which we take as an estimate of $q_t$. In (3.13) we assume for notational simplicity that the minimum exists, however, we do not require that it is unique. If the minimum is not uniquely defined, we can choose as estimate any functions which achieves the minimum and for this function the theoretical results in Section 4 will hold.

REMARK 3.1. It is enough that $\tilde{q}_{n,t}^{w,p}$ is almost minimizer in the sense that

$$
\begin{aligned}
\frac{1}{n_l} &\sum_{i=1}^{n_l} |\tilde{q}_{n,t}^{w,p}(X_{i,t}) - \hat{Y}_{i,t}^{w,\text{new}}|^2 \\
&\leq \min_{h \in \mathcal{H}_{n,p}} \frac{1}{n_l} \sum_{i=1}^{n_l} |h(X_{i,t}) - \hat{Y}_{i,t}^{w,\text{new}}|^2 + o(n^{-1}).
\end{aligned}
\tag{3.14}
$$

The result follows from the proofs.

Let

$$
T_L z = \max\{-L, \min\{L, z\}\}, \qquad z \in \mathbb{R},
\tag{3.15}
$$

denote the truncation operator at threshold level $L > 0$. For a suitable threshold parameter $\beta_n > 0$, to be determined later, we set

$$
\hat{q}_{n,t}^{w,p}(x) = T_{\beta_n} \tilde{q}_{n,t}^{w,p}(x) \qquad (x \in \mathbb{R}^d),
\tag{3.16}
$$

such that $\hat{q}_{n,t}^{w,p}$ is bounded in absolute value by $\beta_n$. Next, we apply the method of splitting the sample to select the parameter $p$; see, for instance, Chapter 7 in [14]. We set

$$
\hat{q}_{n,t}^w(x) = \hat{q}_{n,t}^{w,\hat{p}_t^w}(x) \qquad (x \in \mathbb{R}^d),
\tag{3.17}
$$

where $\hat{p}_t^w \in \mathcal{P}_n$ satisfies

$$
\begin{aligned}
\frac{1}{n_t} &\sum_{i=n_l+1}^{n_l+n_t} |\hat{q}_{n,t}^{w,\hat{p}_t^w}(X_{i,t}) - \hat{Y}_{i,t}^{w,\text{new}}|^2 \\
&= \min_{p \in \mathcal{P}_n} \frac{1}{n_t} \sum_{i=n_l+1}^{n_l+n_t} |\hat{q}_{n,t}^{w,p}(X_{i,t}) - \hat{Y}_{i,t}^{w,\text{new}}|^2.
\end{aligned}
\tag{3.18}
$$

Finally, we explain our choice of $w$. For each $w \in \{0, 1, \ldots, w_{\max}(t)\}$, definition (3.17) provides an estimate $\hat{q}_{n,t}^w$ of $q_t$. The idea is to compute from $\hat{q}_{n,t}^w$ an approximately optimal stopping rule which gives a lower bound on the solution of the optimal stopping problem at time $t$. The optimal candidate for $w$ is the one that maximizes the lower bound. We therefore set

$$
\hat{w}_t = \underset{w \in \{0,1,\ldots,w_{\max}(t)\}}{\arg\max} \frac{1}{n_v} \sum_{i=n_l+n_t+1}^{n} f_{\hat{\tau}_t^w(X_{i,t:T-t-1}^{t,\text{new}})}(X_{i,\hat{\tau}_t^w(X_{i,t:T-t-1}^{t,\text{new}})}^{t,\text{new}}),
\tag{3.19}
$$

where for $w \in \{0, 1, \ldots, w_{\max}(t)\}$ the approximately optimal stopping rule $\hat{\tau}_t^w$ is defined by

$$
\hat{\tau}_t^w = \tau_t(\hat{q}_{n,t}^w, \hat{q}_{n,t+1}, \ldots, \hat{q}_{n,T-2}, \hat{q}_{T-1}),
\tag{3.20}
$$

with $\tau_t(h)$ recursively defined as in (2.10). The specification (3.19) for $\hat{w}_t$ completes the definition of the estimator (3.12).

REMARK 3.2. The quality and the computational cost of the estimator primarily depends on the size of $n_l$, which is used in (3.13) to perform the key nonparametric regression. On the other hand, the magnitude of $n_t$ and $n_v$ is less critical because they are only used to select optimal parameter values from a relatively small discrete set and the corresponding objective functions converge, according to Hoeffding's and Bernstein's inequality, very fast. The impact of $n_t$ and $n_v$ on the overall computation cost is also minor. In practical applications, $n_l$ should be increased as large as affordable by the available computation capacity.

REMARK 3.3. Note that the optimization in (3.18) and (3.19) is performed over a finite set, which implies the existence of an optimizer.

**4. Main theoretical results.** If the stochastic process of the underlying state variables $X_t$ is unbounded, we first localize it to a bounded set $[-A, A]^d$. For many industry models, the localization error can be estimated explicitly. For illustration, we consider a discretely sampled jump-diffusion process $X_t$. Let

$$
\begin{aligned}
\mathbf{G}f(t,x) = (\tfrac{1}{2}\operatorname{tr}(\mathbf{A}\,\nabla^2 f) &+ \langle \mathbf{b}, \nabla f\rangle)(t,x) \\
&+ \int_{\mathbb{R}^d \setminus \{0\}} (f(x+u) - f(x) \\
&\qquad - 1_{\{\|u\|<1\}} \langle u, \nabla f(t,x)\rangle) S(t,x,du)
\end{aligned}
\tag{4.1}
$$

be the generator of the corresponding continuous time process $X_t^0$, where we assume that $\mathbf{A}$, $\mathbf{b}$ are Borel measurable, $\mathbf{A}$ is positive definite, with norms $\|\mathbf{A}\| \leq a_0$, $\|\mathbf{b}\| \leq b_0$, and $S$ is a positive kernel on $\mathbb{R}^d \setminus \{0\}$, Borel measurable in $x$ such that

$$
\sup_x S(t, \|u\|^2 1_{\{\|u\|\leq 1\}} + \|u\| 1_{\{\|u\|>1\}}, du) \leq c_0.
\tag{4.2}
$$

Define

$$
m_t = \sup_{0\leq s\leq t} \|X_s^0 - x\|.
\tag{4.3}
$$

Then, Lemma 17 of [20] states that, for every $\lambda \in \mathbb{R}$ and positive $A$, $\eta$, there exists a constant $k$ only depending on $a_0$ and $c_0$ such that

$$
\mathbf{P}(m_t > A) \leq 2d \exp\left(-\frac{\lambda}{d}(A - \|x\| - b_0 t - \eta) + \frac{\lambda^2}{2} kt(1 + e^{|\lambda|})\right) + \frac{c_0 t}{\eta}.
\tag{4.4}
$$

To localize the process $X_t^0$ to a bounded set $[-A, A]^d$, we replace $X_t^0$ with the process $X_t^{0,A}$ killed at first exit from $[-A, A]^d$. The semi-group of the killed process is

$$
P_t^{0,A} f(x) = \mathbf{E}_x\{f(X_t^{0,A})\} = \mathbf{E}_x\{f(X_t^0) M_t\},
\tag{4.5}
$$

where $M_t$ is the multiplicative functional $M_t = 1_{\{t < \tau_A\}}$ for $\tau_A = \inf\{s \geq 0 | X_s^0 \notin [-A, A]^d\}$; see, for instance, [1]. We obtain

$$
\begin{aligned}
(4.6) \qquad & \sup_{\tau \in \mathcal{T}_{[0,T]}} |\mathbf{E}\{T_L f(X_\tau^0)\} - \mathbf{E}\{T_L f(X_\tau^{0,A})\}| \\
& \leq \sup_{\tau \in \mathcal{T}_{[0,T]}} \mathbf{E}\{T_L f(X_\tau^0) 1_{\{m_\tau > A\}}\} \leq L \mathbf{P}(m_T > A),
\end{aligned}
$$

which, because of (4.4), can be made arbitrarily small by first choosing $\eta$ and then $A$ large enough. Proposition 5.2 in [10] estimates the error if the payoff $f_t$ is replaced by the truncated payoff $T_L f_t$. We arrive at an a priori bound for the localization and payoff truncation error.

In the following we derive the consistency of our estimator (3.12) under the assumption

$$
(4.7) \qquad X_t \in [-A, A]^d \qquad \text{a.s. } (t \in \{0, 1, \ldots, T\}).
$$

In addition, we assume that the payoff $f_s$ is bounded on $[-A, A]^d$ by some constant $L > 0$ such that

$$
(4.8) \qquad |f_s(x)| \leq L \qquad \text{for } x \in [-A, A]^d \text{ and } s \in \{0, 1, \ldots, T\}.
$$

Observe that (4.8) implies $|q_t(x)| \leq L$ for $x \in [-A, A]^d$ and $t \in \{0, 1, \ldots, T\}$, so that $\beta_n = L$ can serve as the truncation parameter for the estimator.

In the sequel we use polynomial splines to define the function spaces $\mathcal{H}_{n,p} = \mathcal{H}_p$ independent of the sample size $n$ and parameterized by $p = (M, \alpha) \in \mathbb{N}_0 \times (0, \infty)$. We note that our results can be extended to other function spaces in a straightforward manner.

For $p = (M, \alpha)$ and $k \in \mathbb{Z}$, we set $u_k = k \cdot \alpha$. Let $B_{k,M} : \mathbb{R} \to \mathbb{R}$ be the univariate B-spline of degree $M$ with knot sequence $(u_l)_{l \in \mathbb{Z}}$ and support $\operatorname{supp}(B_{k,M}) = [u_k, u_{k+M+1}]$. In the case of $M = 0$ the B-spline $B_{k,0}$ is the indicator function of the interval $[u_k, u_{k+1}]$. If $M = 1$, we obtain the so-called hat-functions

$$
B_{k,1}(x) = \begin{cases} \dfrac{x - u_k}{u_{k+1} - u_k}, & \text{for } u_k \leq x \leq u_{k+1}, \\ \dfrac{u_{k+2} - x}{u_{k+2} - u_{k+1}}, & \text{for } u_{k+1} < x \leq u_{k+2}, \\ 0, & \text{else.} \end{cases}
$$

The general definition of $B_{k,M}$ can be found, for example, in [8] or in Section 14.1 of [14]. The B-splines $B_{k,M}$ are basis functions which are piecewise univariate polynomials of degree $M$. They are globally $(M - 1)$-times continuously differentiable, and the $M$th derivative can only jump at the knots $u_l$.

For every multi-index $\mathbf{k} = (k_1, \ldots, k_d) \in \mathbb{Z}^d$, we define the tensor product B-spline $B_{\mathbf{k}, M} : \mathbb{R}^d \to \mathbb{R}$ by

$$B_{\mathbf{k}, M}(x^{(1)}, \ldots, x^{(d)}) = B_{k_1, M}(x^{(1)}) \cdots B_{k_d, M}(x^{(d)}) \qquad (x^{(1)}, \ldots, x^{(d)} \in \mathbb{R}).$$

Let

$$\mathcal{H}_{n,p} = \left\{ \sum_{\mathbf{k} \in \mathbb{Z}^d \, : \, \mathrm{supp}(B_{\mathbf{k}, M}) \cap [-A, A]^d \neq \varnothing} a_{\mathbf{k}} \cdot B_{\mathbf{k}, M} : a_{\mathbf{k}} \in \mathbb{R} \right\}$$

be the span of tensor product B-splines $B_{\mathbf{k}, M}$, such that $\mathrm{supp}(B_{\mathbf{k}, M})$ has a nonempty intersection with $[-A, A]^d$. The spanning functions $B_{\mathbf{k}, M}$ are $(M-1)$-times continuously differentiable, piecewise multivariate polynomial of degree less than or equal to $M$, defined on rectangular domains

$$(4.9) \qquad [u_{k_1}, u_{k_1+1}) \times \cdots \times [u_{k_d}, u_{k_d+1}) \qquad (\mathbf{k} = (k_1, \ldots, k_d) \in \mathbb{Z}^d),$$

and vanish on all of the rectangles (4.9) for which there exists $j \in \{1, \ldots, d\}$ such that either

$$k_j > 0 \quad \text{and} \quad u_{k_j - M} > A$$

or

$$k_j < 0 \quad \text{and} \quad u_{k_j + M + 1} < -A.$$

Consequently, $\mathcal{H}_{n,p}$ is a linear space of functions consisting of piecewise polynomials with respect to equidistant partitions of $\mathbb{R}^d$ into cubes of edge length $\alpha$, vanishing outside a compact set.

For a sample size $n$, we use the parameters

$$\mathcal{P}_n = \{(M, \alpha) : M \in \mathbb{N}_0, M \leq \lceil \log(n) \rceil, \alpha = 2^k \text{ for some } k \in \mathbb{Z}, |k| \leq \lceil \log(n) \rceil \}.$$

Here log denotes the natural logarithm, and for $z \in \mathbb{R}$, we denote by $\lceil z \rceil$ the smallest integer greater than or equal to $z$.

Let $\hat{q}_{n,t}$ be defined as in Section 3 with $\mathcal{P}_n$ and $\mathcal{H}_{n,p}$ as above. Note that $\mathcal{H}_{n,p}$ is a linear function space which implies that the minimum in (3.13) always exists. According to Remark 3.2, the computational cost of the estimator is not adversely affected by large values for $n_t$ and $n_v$ of roughly the size of $n_l$. Therefore, we choose for simplicity $n_v = n_t = \lfloor n/3 \rfloor$ and $n_l = n - n_v - n_t$. Our first result concerns consistency of the estimator.

THEOREM 4.1. *Assume* (4.7) *and* (4.8), *and let the estimate* $\hat{q}_{n,t}$ *be defined as above with* $\beta_n = L$. *Then*

$$\mathbf{E} \int |\hat{q}_{n,t}(x) - q_t(x)|^2 \mu_t(dx) \to 0 \qquad (n \to \infty)$$

*for all* $t \in \{0, 1, \ldots, T-1\}$.

REMARK 4.2. Because convergence in $L_1$ implies convergence in probability, Theorem 4.1 proves, in particular, that $\int |\hat{q}_{n,t} - q_t|^2 \mu_t(dx) \to 0$ in probability as $n \to \infty$.

Next we study the rate of convergence. It is well known in nonparametric regression that without smoothness assumptions on the regression function the rate of convergence can be arbitrarily slow (cf., e.g., [7, 9] or [14], Chapter 3). We assume that the continuation values $q_t$ are $(p, C)$-smooth according to the following definition.

DEFINITION 4.3. Let $p = k + \beta$ for some $k \in \mathbb{N}_0$, $\beta \in (0, 1]$, and let $C > 0$. A function $f : \mathbb{R}^d \to \mathbb{R}$ is called $(p, C)$-smooth, if all partial derivatives

$$\frac{\partial f}{\partial^{\alpha_1} x^{(1)} \ldots \partial^{\alpha_d} x^{(d)}}$$

of total order $\alpha_1 + \cdots + \alpha_d = k$ exist and satisfy

$$\left| \frac{\partial f}{\partial^{\alpha_1} x^{(1)} \ldots \partial^{\alpha_d} x^{(d)}}(x) - \frac{\partial f}{\partial^{\alpha_1} x^{(1)} \ldots \partial^{\alpha_d} x^{(d)}}(z) \right| \leq C \cdot \|x - z\|^{\beta}$$

for all $x, z \in \mathbb{R}^d$.

Such a smoothness assumption is not unreasonable. For a sufficiently regular diffusion or jump-diffusion process, the semi-group of Markov transition operators $P_{s,t}(g)(x) = \mathbf{E}[g(X_t)|X_s = x]$ is strongly smoothing already for arbitrarily small time steps. In particular, we can expect that the continuation value $q_t = P_{t,t+1}(\max((f_{t+1}, q_{t+1}))$ is $(p, C)$-smooth under suitable assumptions on $X_t$ and the payoff $f_t$. At this point, it also becomes clear why it is unfavorable to directly work with the value function $v_t$ which does not retain the smoothness because the maximum operation is applied after the transition operator. Next, we address the rate of convergence of the estimator.

THEOREM 4.4. Let $p = k + \beta$ for some $k \in \mathbb{N}_0$, $\beta \in (0, 1]$, and let $C > 0$. Assume $k \leq M_{\max}$, (4.7), (4.8) and

$$q_t \quad (p, C)\text{-smooth}$$

for all $t \in \{0, 1, \ldots, T - 1\}$. Let the estimate $\hat{q}_{n,t}$ be defined as above with $\beta_n = L$. Then for every $t \in \{0, 1, \ldots, T - 1\}$,

$$\mathbf{E} \int |\hat{q}_{n,t}(x) - q_t(x)|^2 \mu_t(dx) \leq const \cdot C^{2d/(2p+d)} \cdot \left( \frac{\log n}{n} \right)^{2p/(2p+d)}.$$

REMARK 4.5. We would like to stress that in Theorems 4.1 and 4.4 there is no assumption on the distribution of $X$ besides the assumption (4.7). In particular, it is not required that $X_t$ has a density with respect to the Lebesgue–Borel measure.

REMARK 4.6. It is well known that the optimal rate of convergence for the estimation of $(p, C)$-smooth functions is $n^{-2p/(2p+d)}$ (see, e.g., [27] or [14], Chapter 3). Hence, the rate of convergence in Theorem 4.4 is optimal up to a logarithmic factor.

REMARK 4.7. The definition of the estimator in Theorem 4.4 does not depend on the degree of smoothness of $q_t$ represented by $(p, C)$. Nevertheless, the estimator achieves the optimal rate of convergence for a particular smoothness of the continuation value. In this sense the estimator is able to adapt automatically to the smoothness of the continuation value, in contrast to the estimates in [10].

REMARK 4.8. Assume $X_0 = x_0$ a.s. for some $x_0 \in [-A, A]^d$. We can estimate the price

$$V_0 = v_0(x_0) = \max\{f_0(x_0), q_0(x_0)\}$$

[cf. (1.2), (2.2) and (2.5)] of the Bermudan option by

$$\hat{V}_0 = \max\{f_0(x_0), \hat{q}_{n,0}(x_0)\}.$$

Since the distribution of $X_0$ is concentrated at $x_0$, Theorem 4.4 leads to the error bound

$$
\begin{aligned}
\mathbf{E}\{|\hat{V}_0 - V_0|^2\} &= \mathbf{E}\{|\max\{f_0(x_0), \hat{q}_{n,0}(x_0)\} - \max\{f_0(x_0), q_0(x_0)\}|^2\} \\
&\leq \mathbf{E}\{|\hat{q}_{n,0}(x_0) - q_0(x_0)|^2\} \\
&\leq const \cdot C^{2d/(2p+d)} \cdot \left(\frac{\log n}{n}\right)^{2p/(2p+d)}.
\end{aligned}
$$

**5. Finite sample behavior.** In this section we illustrate the finite sample behavior of our algorithm (EKT) in comparison to the Longstaff–Schwartz (LS) and Tsitsiklis–Van Roy (TR) algorithm. To compare the three algorithms, we proceed as follows. We independently generate sample paths and compute for each algorithm the Monte Carlo estimates (MCE) of the price (1.6). Because all three algorithms provide a lower bound for the optimal stopping value, and since we evaluate the approximative optimal stopping rule on independent sets of sample path, a higher MCE indicates a better performance.

The underlying model for the dynamics of the stocks is a simple geometric Brownian motion. We apply a Euler scheme to discretize the time interval $[0, 1]$ into $m$ time steps. Consequently, the prices of the underlying stocks on the time grid $0, \frac{1}{m}, \ldots, \frac{m-1}{m}, 1$ are given by

$$X_{i,j} = X_0 \cdot \exp\left( \left( r - \frac{1}{2}\sigma^2 \right) \cdot \frac{j}{m} + \frac{\sigma}{\sqrt{m}} \cdot W_{i,j} \right)$$

(5.1)
$$(i = 1, \ldots, n, j = 1, \ldots, m).$$

Here, $X_0$ is the initial stock price at time 0, $r$ is the risk-free interest rate, $\sigma$ the instantaneous volatility, and

$$W_{i,j} = \sum_{l=1}^{j} Z_{i,l}$$

is the sum of independent standard normally distributed random variables $Z_{i,l} (i = 1, \ldots, n, l = 1, \ldots, m)$. All option contracts are based on a time to maturity of 1 year and a risk-free continuously compounded interest rate $r = 0.05$.

Figures 1 and 3 report the results for 100 independent MCE of ordinary Bermudan put option and for a more complicated Bermudan option with a strangle spread payoff. Each algorithm is based on a sample size $n = 10000$. For (LS) and (TR), we use polynomials of degree 3. For (EKT), we set the number of learning, training and validation samples to $n_l = 6000$, $n_t = 2000$ and $n_v = 2000$, and choose the degree $M$, the knot distance $\alpha$ and the look-ahead parameter $w(t)$ in a data-dependent manner as described in Section 3 from the sets $M \in \{0, 1, 2\}$, $\alpha \in \{\frac{100}{2}, \frac{100}{2^2}, \frac{100}{2^3}, \frac{100}{2^4}\}$, and $w(t) \in \{0, 4, T - t - 1\}$.

We first analyze the results in Figure 1 for a Bermudan put with exercise price 90 on an underlying with instantaneous volatility $\sigma = 0.25$. The time discretization is performed in monthly steps. Our algorithm is slightly better than (LS) and comparable to (TR). This is not surprising, since it is well known that for simple payoff functions both (LS) and (TR) perform rather very well.

Figure 3 consolidates the simulation results of a Bermudan option with strangle spread payoff with 50, 90, 110 and 150, as illustrated in Figure 2. The volatility is increased to $\sigma = 0.5$, the time discretization is set to $m = 48$. This time (EKT) provides a higher MCE of the option price and therefore clearly outperforms (LS) and (TR).

Finally, Figure 4 reports the simulation results of a Bermudan basket option with strangle spread payoff on the average of three correlated underlyings. The option prices are normalized to start at 1. The strikes are set at 0.85, 0.95, 1.05 and 1.15. This time (EKT) is based on degrees $M \in \{0, 1, 2\}$, knot distance $\alpha \in \{1, 1.5, 2, 4\}$ and a reduced sample size of only $n = 4000$, split into $n_t = 800$, $n_l = 2400$ and $n_v = 800$. (LS) and (TR) still use $n = 10000$ but approximate the continuation value with polynomials of degree 2 (as polynomials of degree 3 resulted in lower MCE). Again, (EKT) provides the highest MCE of the option price.

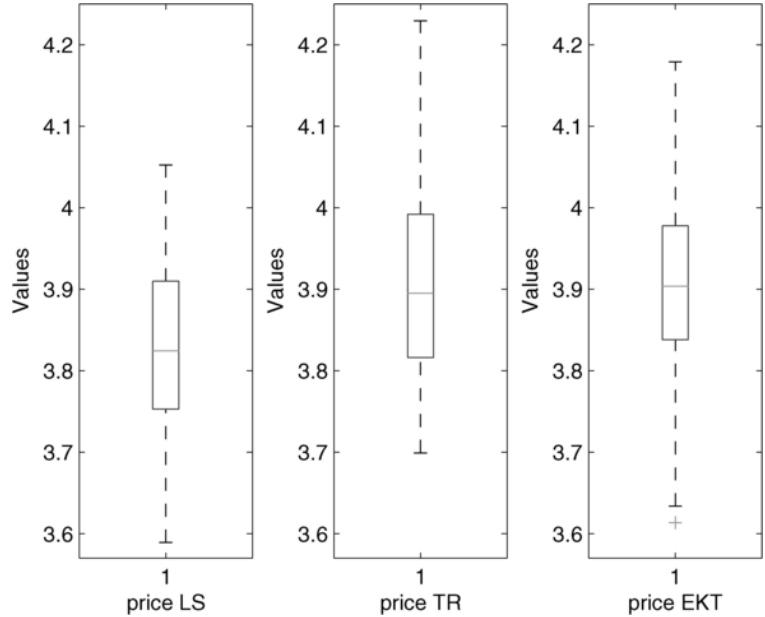

Fig. 1. *Realized option prices of Bermudan put option. The boxes stretch from the 25th percentile to the 75th percentile, the median is shown as a line across the box.*

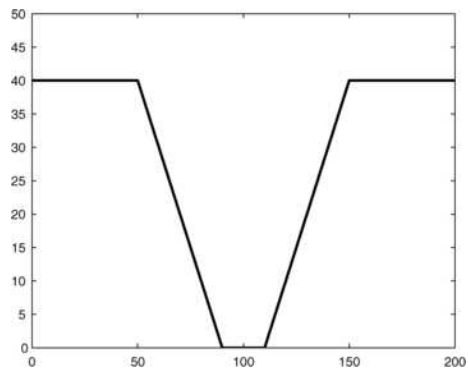

Fig. 2. *Strangle spread payoff with strike prices 50, 90, 110 and 150.*

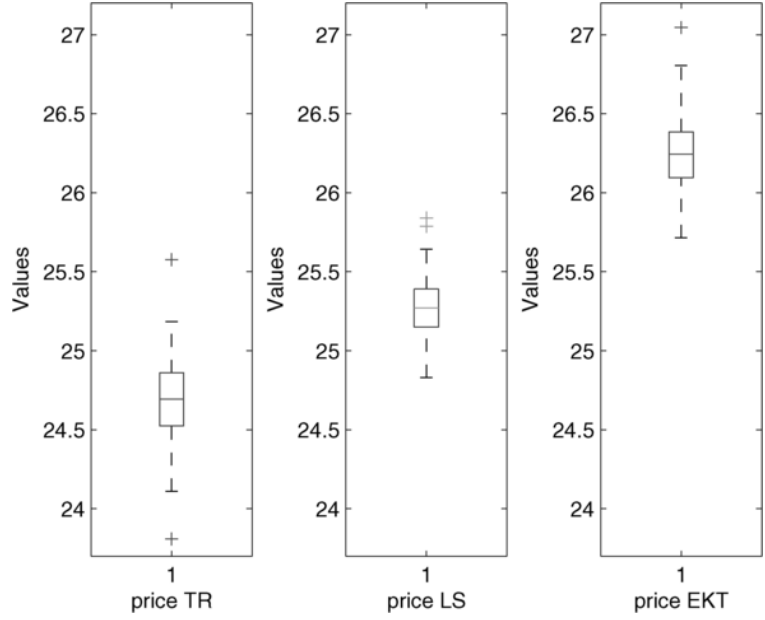

Fig. 3. *Realized option prices of Bermudan option with strangle spread-payoff.*

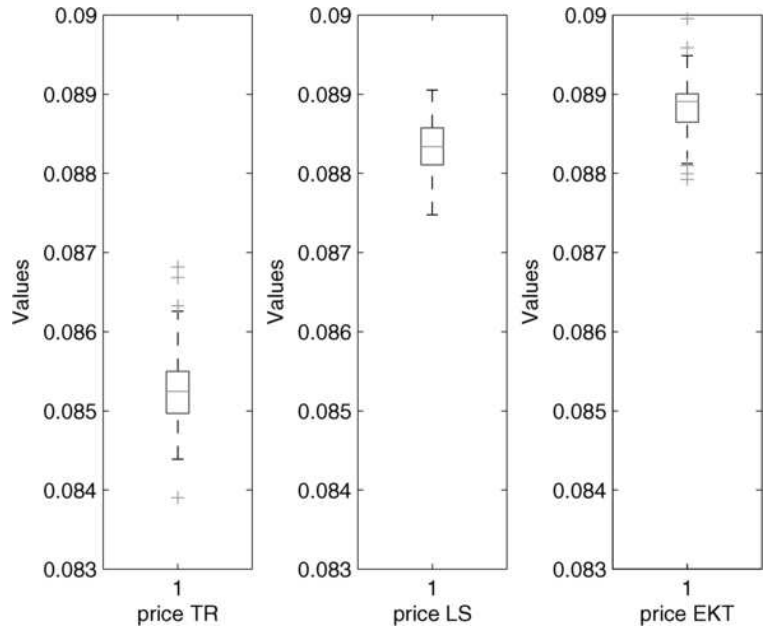

Fig. 4. *Realized option prices of Bermudan basket option with strangle spread-payoff based on the average of three correlated underlyings.*

**6. Proofs.** In the proofs we will need an auxiliary result on the properties of the method of splitting the sample, which we formulate and prove for the sake of generality in a fixed design regression model.

Let $x_1, \ldots, x_n \in \mathbb{R}^d$ and let $Y_1, \ldots, Y_n$ be independent square integrable random variables which satisfy

$$\mathbf{E} Y_i = m(x_i) \qquad (i = 1, \ldots, n)$$

for some function $m : \mathbb{R}^d \to \mathbb{R}$. Let $\mathcal{P}_n$ be a finite set of parameters and assume that for each $p \in \mathcal{P}_n$ an estimate $m_p : \mathbb{R}^d \to \mathbb{R}$ is given. Choose $p^* \in \mathcal{P}_n$ by minimizing the empirical $L_2$ risk on the sample $(x_1, Y_1), \ldots, (x_n, Y_n)$, that is, assume

$$\frac{1}{n} \sum_{i=1}^{n} |m_{p^*}(x_i) - Y_i|^2 = \min_{p \in \mathcal{P}_n} \frac{1}{n} \sum_{i=1}^{n} |m_p(x_i) - Y_i|^2.$$

Then, the following bound on the error

$$\frac{1}{n} \sum_{i=1}^{n} |m_{p^*}(x_i) - m(x_i)|^2$$

of $m_{p^*}$ holds.

LEMMA 6.1. *Under the above assumptions, we have for each $\epsilon > 0$*

$$\mathbf{P} \left\{ \frac{1}{n} \sum_{i=1}^{n} |m_{p^*}(x_i) - m(x_i)|^2 > \epsilon + 18 \cdot \min_{p \in \mathcal{P}_n} \frac{1}{n} \sum_{i=1}^{n} |m_p(x_i) - m(x_i)|^2 \right\}$$

$$\leq c_1 \cdot \max_{i=1, \ldots, n} \mathbf{E} Y_i^2 \cdot \frac{|\mathcal{P}_n|}{\epsilon \cdot n}$$

*for some constant $c_1$ which does not depend on $n$ or $\epsilon$.*

PROOF. Set

$$m^* = \arg \min_{f \in \{m_p : p \in \mathcal{P}_n\}} \frac{1}{n} \sum_{i=1}^{n} |f(x_i) - m(x_i)|^2.$$

By Lemma 1 in [17] or standard results from the book [30] (see proof of Theorem 10.11 in [30]), we have

$$\mathbf{P} \left\{ \frac{1}{n} \sum_{i=1}^{n} |m_{p^*}(x_i) - m(x_i)|^2 > \epsilon + 18 \cdot \min_{p \in \mathcal{P}_n} \frac{1}{n} \sum_{i=1}^{n} |m_p(x_i) - m(x_i)|^2 \right\}$$

$$\leq \mathbf{P} \left\{ \frac{\epsilon}{2} < \frac{1}{n} \sum_{i=1}^{n} |m_{p^*}(x_i) - m^*(x_i)|^2 \right.$$

$$\leq \frac{16}{n} \sum_{i=1}^{n} (m_{p^*}(x_i) - m^*(x_i)) \cdot (Y_i - m(x_i)) \bigg\}$$

$$\leq |\mathcal{P}_n| \cdot \max_{p \in \mathcal{P}_n} \mathbf{P} \bigg\{ \frac{\epsilon}{2} < \frac{1}{n} \sum_{i=1}^{n} |m_p(x_i) - m^*(x_i)|^2$$

$$\leq \frac{16}{n} \sum_{i=1}^{n} (m_p(x_i) - m^*(x_i)) \cdot (Y_i - m(x_i)) \bigg\}$$

$$\leq |\mathcal{P}_n| \cdot \max_{p \in \mathcal{P}_n} \sum_{s=0}^{\infty} \mathbf{P} \bigg\{ 2^{s-1} \epsilon < \frac{1}{n} \sum_{i=1}^{n} |m_p(x_i) - m^*(x_i)|^2 \leq 2^s \epsilon,$$

$$\frac{1}{n} \sum_{i=1}^{n} |m_p(x_i) - m^*(x_i)|^2$$

$$\leq \frac{16}{n} \sum_{i=1}^{n} (m_p(x_i) - m^*(x_i)) \cdot (Y_i - m(x_i)) \bigg\}$$

$$\leq |\mathcal{P}_n| \cdot \sum_{s=0}^{\infty} \max_{\substack{p \in \mathcal{P}_n \\ (1/n) \sum_{i=1}^{n} |m_p(x_i)-m^*(x_i)|^2 \leq 2^s \epsilon}} \mathbf{P} \bigg\{ \frac{1}{n} \sum_{i=1}^{n} (m_p(x_i) - m^*(x_i))$$

$$\cdot (Y_i - m(x_i)) > \frac{2^s \epsilon}{32} \bigg\}.$$

Because of the variance estimate

$$\mathbf{V} \bigg( \frac{1}{n} \sum_{i=1}^{n} (m_p(x_i) - m^*(x_i)) \cdot (Y_i - m(x_i)) \bigg)$$

$$\leq \frac{1}{n^2} \sum_{i=1}^{n} (m_p(x_i) - m^*(x_i))^2 \cdot \max_{i=1,\dots,n} \mathbf{E} Y_i^2,$$

we can bound the right-hand side from above with Chebyshev's inequality by

$$|\mathcal{P}_n| \cdot \sum_{s=0}^{\infty} \frac{(1/n) \cdot 2^s \cdot \epsilon \cdot \max_{i=1,\dots,n} \mathbf{E} Y_i^2}{(2^s \epsilon/32)^2} = \frac{|\mathcal{P}_n|}{n} \cdot \frac{\max_{i=1,\dots,n} \mathbf{E} Y_i^2}{\epsilon} \cdot \sum_{s=0}^{\infty} \frac{32^2}{2^s}. \quad \Box$$

PROOF OF THEOREM 4.1.  Because of

$$\mathbf{E} \int |\hat{q}_{n,t}(x) - q_t(x)|^2 \mu_t(dx) \leq \sum_{w=0}^{w_{\max}(t)} \mathbf{E} \int |\hat{q}_{n,t}^w(x) - q_t(x)|^2 \mu_t(dx),$$

it is enough to prove that

$$(6.1) \qquad \mathbf{E} \int |\hat{q}_{n,t}^w(x) - q_t(x)|^2 \mu_t(dx) \to 0 \qquad (n \to \infty)$$

for every $t \in \{0, 1, \ldots, T-1\}$ and every $w \in \{0, 1, \ldots, w_{\max}(t)\}$.

Fix $t \in \{0, 1, \ldots, T-1\}$ and assume (by induction) that we have for every $s \in \{t+1, \ldots, T-1\}$ and every $v \in \{0, 1, \ldots, w_{\max}(s)\}$

$$(6.2) \qquad \mathbf{E} \int |\hat{q}_{n,s}^v(x) - q_s(x)|^2 \mu_t(dx) \to 0 \qquad (n \to \infty).$$

Fix $w \in \{0, 1, \ldots, w_{\max}(t)\}$. In the following we show

$$(6.3) \qquad \mathbf{E} \int |\hat{q}_{n,t}^w(x) - q_t(x)|^2 \mu_t(dx) \to 0 \qquad (n \to \infty).$$

To this end, we apply for a fixed $p_n \in \mathcal{P}_n$ the error decomposition

$$\int |\hat{q}_{n,t}^w(x) - q_t(x)|^2 \mu_t(dx)$$

$$= \int |\hat{q}_{n,t}^w(x) - q_t(x)|^2 \mu_t(dx) - \frac{1}{n_t} \sum_{i=n_l+1}^{n_l+n_t} |\hat{q}_{n,t}^w(X_{i,t}) - q_t(X_{i,t})|^2$$

$$+ \left( \frac{1}{n_t} \sum_{i=n_l+1}^{n_l+n_t} |\hat{q}_{n,t}^w(X_{i,t}) - q_t(X_{i,t})|^2 \right.$$

$$\left. - \frac{2}{n_t} \sum_{i=n_l+1}^{n_l+n_t} |\hat{q}_{n,t}^w(X_{i,t}) - q_t^{w,\text{new}}(X_{i,t})|^2 \right)$$

$$+ \left( \frac{2}{n_t} \sum_{i=n_l+1}^{n_l+n_t} |\hat{q}_{n,t}^w(X_{i,t}) - q_t^{w,\text{new}}(X_{i,t})|^2 \right.$$

$$\left. - \frac{36}{n_t} \sum_{i=n_l+1}^{n_l+n_t} |\hat{q}_{n,t}^{w,p_n}(X_{i,t}) - q_t^{w,\text{new}}(X_{i,t})|^2 \right)$$

$$+ \left( \frac{36}{n_t} \sum_{i=n_l+1}^{n_l+n_t} |\hat{q}_{n,t}^{w,p_n}(X_{i,t}) - q_t^{w,\text{new}}(X_{i,t})|^2 \right.$$

$$\left. - \frac{72}{n_t} \sum_{i=n_l+1}^{n_l+n_t} |\hat{q}_{n,t}^{w,p_n}(X_{i,t}) - q_t(X_{i,t})|^2 \right)$$

$$+ \frac{72}{n_t} \sum_{i=n_l+1}^{n_l+n_t} |\hat{q}_{n,t}^{w,p_n}(X_{i,t}) - q_t(X_{i,t})|^2$$

$$= \sum_{j=1}^{5} T_{j,n}.$$

The proof will be completed once we have shown that

$$(6.4) \qquad \limsup_{n \to \infty} \mathbf{E} T_{j,n} \le 0$$

for $j \in \{1, 2, \ldots, 5\}$.

From now on we denote by $\mathcal{D}_{n,t+1}^T$ the set of all the data used in the construction of the estimates $\hat{q}_{n,s}^{w,p}$ for $s > t$, $w \in \{0, 1, \ldots, w_{\max}(s)\}$ and $p \in \mathcal{P}_n$.

Because $\hat{q}_{n,t}^w$ and $q_t$ are bounded in absolute value by $L$, we conclude from Hoeffding's inequality (see, e.g., Lemma A.3 in [14]) that

$$\mathbf{P}\{T_{1,n} > \epsilon | X_{i,t:w_{\max}(t)+1}^{t,\text{new}} \ (i = 1, \ldots, n_l), \mathcal{D}_{n,t+1}\}$$

$$\leq |\mathcal{P}_n| \cdot \max_{p \in \mathcal{P}_n} \mathbf{P}\bigg\{ \int |\hat{q}_{n,t}^{w,p}(x) - q_t(x)|^2 \mu_t(dx)$$

$$- \frac{1}{n_t} \sum_{i=n_l+1}^{n_l+n_t} |\hat{q}_{n,t}^{w,p}(X_{i,t}) - q_t(X_{i,t})|^2 > \epsilon$$

$$\bigg| X_{i,t:w_{\max}(t)+1}^{t,\text{new}} (i = 1, \ldots, n_l), \mathcal{D}_{n,t+1}^T \bigg\}$$

$$\leq |\mathcal{P}_n| \cdot \exp\bigg( -\frac{2n_t\epsilon^2}{(4L^2)^2} \bigg) = \exp\bigg( \log(|\mathcal{P}_n|) - \frac{2n_t\epsilon^2}{16L^4} \bigg).$$

Thus,

$$\mathbf{E}T_{1,n} \leq \int_0^\infty \mathbf{P}\{T_{1,n} > s\} \, ds$$

$$= \int_0^\infty \mathbf{E}\{\mathbf{P}\{T_{1,n} > s | X_{i,t:w_{\max}(t)+1}^{t,\text{new}}(i = 1, \ldots, n_l), \mathcal{D}_{n,t+1}^T\}\} \, ds$$

$$\leq 4L^2\sqrt{\log(|\mathcal{P}_n|)/n_t} + \int_{4L^2\sqrt{\log(|\mathcal{P}_n|)/n_t}}^\infty \exp\bigg( -\frac{n_t s^2}{16L^4} \bigg) \, ds$$

$$\leq 4L^2\sqrt{\log(|\mathcal{P}_n|)/n_t}$$

$$+ \int_{4L^2\sqrt{\log(|\mathcal{P}_n|)/n_t}}^\infty \exp\bigg( -\frac{n_t \cdot 4L^2\sqrt{\log(|\mathcal{P}_n|)/n_t}}{16L^2} \cdot s \bigg) \, ds$$

$$\leq 4L^2\sqrt{\log(|\mathcal{P}_n|)/n_t}$$

$$+ \frac{4L^2}{n_t\sqrt{\log(|\mathcal{P}_n|)/n_t}} \cdot \exp\bigg( -\log(|\mathcal{P}_n|) \bigg) \to 0 \qquad (n \to \infty).$$

Furthermore, by $a^2 = (a - b + b)^2 \leq 2(a - b)^2 + 2b^2$, we get

$$T_{2,n} \leq \frac{2}{n_t} \sum_{i=n_l+1}^{n_l+n_t} |q_t^{w,\text{new}}(X_{i,t}) - q_t(X_{i,t})|^2,$$

from which we conclude, together with (3.7) and (6.2), that

$$\mathbf{E}T_{2,n} = \mathbf{E}\{\mathbf{E}\{T_{2,n}|X^{t,\mathrm{new}}_{i,t:w_{\max}(t)+1} \ (i=1,\dots,n_l), \mathcal{D}^T_{n,t+1}\}\}$$

$$\leq 2\mathbf{E}\int |q^{w,\mathrm{new}}_t(x) - q_t(x)|^2 \mu_t(dx) \to 0 \qquad (n \to \infty).$$

In a similar way we obtain

$$\mathbf{E}T_{4,n} \leq 72\mathbf{E}\int |q^{w,\mathrm{new}}_t(x) - q_t(x)|^2 \mu_t(dx) \to 0 \qquad (n \to \infty).$$

To bound $T_{3,n}$, we use Lemma 6.1, which shows

$$\mathbf{P}\{T_{3,n} > \epsilon | X^{t,\mathrm{new}}_{i,t:w_{\max}(t)+1} \ (i=1,\dots,n_l), \mathcal{D}^T_{n,t+1}\}$$

$$\leq \mathbf{P}\Big\{\frac{1}{n_t}\sum_{i=n_l+1}^{n_l+n_t} |\hat{q}^w_{n,t}(X_{i,t}) - q^{w,\mathrm{new}}_t(X_{i,t})|^2$$

$$> \frac{\epsilon}{2} + 18 \cdot \min_{p \in \mathcal{P}_n} \frac{1}{n_t}\sum_{i=n_l+1}^{n_l+n_t} |\hat{q}^{w,p}_{n,t}(X_{i,t}) - q^{w,\mathrm{new}}_t(X_{i,t})|^2$$

$$\Big| X^{t,\mathrm{new}}_{i,t:w_{\max}(t)+1} \ (i=1,\dots,n_l), \mathcal{D}^T_{n,t+1}\Big\}$$

$$\leq c_2 \cdot \frac{|\mathcal{P}_n|}{\epsilon \cdot n_t}.$$

This implies for any $u > 0$ that

$$\mathbf{E}T_{3,n} \leq \int_0^\infty \mathbf{P}\{T_{3,n} > \epsilon\}\, d\epsilon$$

$$\leq \int_0^\infty \mathbf{E}\{\mathbf{P}\{T_{3,n} > \epsilon | X^{t,\mathrm{new}}_{i,t:w_{\max}(t)+1} \ (i=1,\dots,n_l), \mathcal{D}^T_{n,t+1}\}\}\, d\epsilon$$

$$\leq u + \int_u^{const} c_2 \cdot \frac{|\mathcal{P}_n|}{\epsilon \cdot n_t}\, d\epsilon$$

$$= u + c_2 \cdot \frac{|\mathcal{P}_n|}{n_t} \cdot (\log(const) - \log u).$$

To get to the last line, we have used that (3.16) and the boundedness of $q^{w,\mathrm{new}}_t$ (which is a consequence of the boundedness of $f_t$ on $[-A,A]^d$) yield

$$T_{3,n} \leq \frac{2}{n_t}\sum_{i=n_l+1}^{n_l+n_t} |\hat{q}^w_{n,t}(X_{i,t}) - q^{w,\mathrm{new}}_t(X_{i,t})|^2 \leq const.$$

Setting $u = |\mathcal{P}_n|/n_t$, we arrive at

$$\limsup_{n \to \infty} \mathbf{E}T_{3,n} \leq 0.$$

Furthermore,

$$
\begin{aligned}
(6.5) \quad \mathbf{E}T_{5,n} &= \mathbf{E}\{\mathbf{E}\{T_{5,n}|X_{i,t:w_{\max}(t)+1}^{t,\mathrm{new}}\ i=1,\ldots,n_l),\mathcal{D}_{n,t+1}^T\}\} \\
&= 72\cdot\mathbf{E}\int|\hat{q}_{n,t}^{w,p_n}(x)-q_t(x)|^2\mu_t(dx).
\end{aligned}
$$

Consequently, it remains to verify that

$$
(6.6) \qquad \mathbf{E}\int|\hat{q}_{n,t}^{w,p_n}(x)-q_t(x)|^2\mu_t(dx)\to 0 \qquad (n\to\infty)
$$

for some suitably selected $p_n\in\mathcal{P}_n$.

To prove (6.6), we set $p_n=(0,2^{-\lceil\log_2(n)/(2+d)\rceil})$ (where $\log_2$ is the logarithm for base 2) and consider the error decomposition

$$
\begin{aligned}
&\int|\hat{q}_{n,t}^{w,p_n}(x)-q_t(x)|^2\mu_t(dx) \\
&\quad = \int|\hat{q}_{n,t}^{w,p_n}(x)-q_t(x)|^2\mu_t(dx)-\frac{2}{n_l}\sum_{i=1}^{n_l}|\hat{q}_{n,t}^{w,p_n}(X_{i,t})-q_t(X_{i,t})|^2 \\
&\qquad +\frac{2}{n_l}\sum_{i=1}^{n_l}|\hat{q}_{n,t}^{w,p_n}(X_{i,t})-q_t(X_{i,t})|^2-\frac{2}{n_l}\sum_{i=1}^{n_l}|\tilde{q}_{n,t}^{w,p_n}(X_{i,t})-q_t(X_{i,t})|^2 \\
&\qquad +\left(\frac{2}{n_l}\sum_{i=1}^{n_l}|\tilde{q}_{n,t}^{w,p_n}(X_{i,t})-q_t(X_{i,t})|^2\right. \\
&\qquad\qquad \left.-\frac{4}{n_l}\sum_{i=1}^{n_l}|\tilde{q}_{n,t}^{w,p_n}(X_{i,t})-q_t^{w,\mathrm{new}}(X_{i,t})|^2\right) \\
&\qquad +\frac{4}{n_l}\sum_{i=1}^{n_l}|\tilde{q}_{n,t}^{w,p_n}(X_{i,t})-q_t^{w,\mathrm{new}}(X_{i,t})|^2 \\
&\quad = \sum_{j=6}^{9}T_{j,n}.
\end{aligned}
$$

Because $q_t$ is bounded in absolute value by $L$, we have

$$
T_{7,n}\le 0 \quad\text{and}\quad \mathbf{E}T_{7,n}\le 0.
$$

In the same way as for $T_{2,n}$, we obtain from (3.7) and (6.2)

$$
\begin{aligned}
\mathbf{E}T_{8,n} &\le 4\cdot\mathbf{E}\left\{\mathbf{E}\left\{\frac{1}{n_l}\sum_{i=1}^{n_l}|q_t(X_{i,t})-q_t^{w,\mathrm{new}}(X_{i,t})|^2\bigg|\mathcal{D}_{n,t+1}^T\right\}\right\} \\
&= 4\cdot\mathbf{E}\int|q_t(x)-q_t^{w,\mathrm{new}}(x)|^2\mu_t(dx)\to 0 \qquad (n\to\infty),
\end{aligned}
$$

where the last equality follows from the fact that the conditional expectation $q_t^{w,\text{new}}(x)$ does not depend on data from time $t$.

Next, we estimate $T_{6,n}$. The functions $\hat{q}_{n,t}^{w,p_n}$ and $q_t$ are bounded in absolute value by $L$, and $\tilde{q}_{n,t}^{w,p_n}$ belongs to the linear vector space $\mathcal{H}_{n,p_n}$, whose dimension $D_n$ is bounded by some constant (depending on $A$) times $n^{d/(2+d)}$. As in the proof of Theorem 11.3 in [14] [in particular, the proof of inequality (11.6)], we obtain

$$\mathbf{E}T_{6,n} = \mathbf{E}\{\mathbf{E}\{T_{6,n}|\mathcal{D}_{n,t+1}^T\}\} \leq c_3 L^2 \frac{(\log n_l + 1) \cdot n^{d/(2+d)}}{n_l} \to 0 \qquad (n \to \infty).$$

It remains to bound $T_{9,n}$. With

$$\sigma^2 = \sup_{x \in \mathbb{R}^d} \mathbf{E}^*\{|\hat{Y}_{1,t}^{w,\text{new}}|^2|X_{1,t} = x\} \leq 4L^2 < \infty,$$

we conclude from Theorem 11.1 in [14]

$$\mathbf{E}\{T_{9,n}|X_{i,t} \ (i=1,\ldots,n_l), \mathcal{D}_{n,t+1}^T\}$$

$$\leq 4\sigma^2 \frac{c_4 n^{d/(2+d)}}{n_l} + 4 \min_{h \in \mathcal{H}_{n,p_n}} \frac{1}{n_l} \sum_{i=1}^{n_l} |h(X_{i,t}) - q_t^{w,\text{new}}(X_{i,t})|^2,$$

which then leads to

$$\mathbf{E}T_{9,n} = \mathbf{E}\{\mathbf{E}\{T_{9,n}|X_{i,t} \ (i=1,\ldots,n_l), \mathcal{D}_{n,t+1}^T\}\}$$

$$\leq 4\sigma^2 \frac{c_4 n^{d/(2+d)}}{n_l} + 4 \min_{h \in \mathcal{H}_{n,p_n}} \mathbf{E}\int |h(x) - q_t^{w,\text{new}}(x)|^2 \mu_t(dx)$$

$$\leq 4\sigma^2 \frac{c_4 n^{d/(2+d)}}{n_l} + 8\mathbf{E}\int |q_t^{w,\text{new}}(x) - q_t(x)|^2 \mu_t(dx)$$

$$+ 8 \min_{h \in \mathcal{H}_{n,p_n}} \int |h(x) - q_t(x)|^2 \mu_t(dx).$$

Because of (3.7), (6.2) and

$$\int |q_t(x)|^2 \mu_t(dx) \leq L^2 < \infty,$$

which implies that $q_t$ can be approximated arbitrarily closely by functions from $\mathcal{H}_{n,p_n}$ (this is a consequence of Theorem A.1 in [14] and the fact that any continuous function can be approximated in the supremum norm on the compact set $[-A, A]^d$ arbitrarily closely by the piecewise constant functions in $\mathcal{H}_{n,p_n}$ as $n \to \infty$), the right-hand side of the above inequality tends to zero for $n \to \infty$. The proof of Theorem 4.1 is complete. $\square$

PROOF OF THEOREM 4.4. The proof is similar to the proof of Theorem 4.1. The main difference is that we use Bernstein's inequality instead of Hoeffding's inequality, which requires that we also control the variance. Because of

$$\mathbf{E} \int |\hat{q}_{n,t}(x) - q_t(x)|^2 \mu_t(dx) \leq \sum_{w=0}^{w_{\max}(t)} \mathbf{E} \int |\hat{q}_{n,t}^w(x) - q_t(x)|^2 \mu_t(dx),$$

it suffices to show

$$(6.7) \quad \mathbf{E} \int |\hat{q}_{n,t}^w(x) - q_t(x)|^2 \mu_t(dx) \leq const \cdot C^{2d/(2p+d)} \cdot \left( \frac{\log n}{n} \right)^{2p/(2p+d)},$$

for every $t \in \{0, 1, \ldots, T-1\}$ and every $w \in \{0, 1, \ldots, w_{\max}(t)\}$.

Fix $t \in \{0, 1, \ldots, T-1\}$ and assume (by induction) that we have for every $s \in \{t+1, \ldots, T-1\}$ and every $v \in \{0, 1, \ldots, w_{\max}(s)\}$

$$(6.8) \quad \mathbf{E} \int |\hat{q}_{n,s}^v(x) - q_s(x)|^2 \mu_t(dx) \leq const \cdot C^{2d/(2p+d)} \cdot \left( \frac{\log n}{n} \right)^{2p/(2p+d)}.$$

Fix $w \in \{0, 1, \ldots, w_{\max}(t)\}$. We show

$$(6.9) \quad \mathbf{E} \int |\hat{q}_{n,t}^w(x) - q_t(x)|^2 \mu_t(dx) \leq const \cdot C^{2d/(2p+d)} \cdot \left( \frac{\log n}{n} \right)^{2p/(2p+d)}.$$

To this end, we apply for fixed $p_n \in \mathcal{P}_n$ the error decomposition

$$\int |\hat{q}_{n,t}^w(x) - q_t(x)|^2 \mu_t(dx)$$

$$= \int |\hat{q}_{n,t}^w(x) - q_t(x)|^2 \mu_t(dx) - \frac{2}{n_t} \sum_{i=n_l+1}^{n_l+n_t} |\hat{q}_{n,t}^w(X_{i,t}) - q_t(X_{i,t})|^2$$

$$+ \left( \frac{2}{n_t} \sum_{i=n_l+1}^{n_l+n_t} |\hat{q}_{n,t}^w(X_{i,t}) - q_t(X_{i,t})|^2 \right.$$

$$\left. - \frac{4}{n_t} \sum_{i=n_l+1}^{n_l+n_t} |\hat{q}_{n,t}^w(X_{i,t}) - q_t^{w,\text{new}}(X_{i,t})|^2 \right)$$

$$+ \left( \frac{4}{n_t} \sum_{i=n_l+1}^{n_l+n_t} |\hat{q}_{n,t}^w(X_{i,t}) - q_t^{w,\text{new}}(X_{i,t})|^2 \right.$$

$$\left. - \frac{72}{n_t} \sum_{i=n_l+1}^{n_l+n_t} |\hat{q}_{n,t}^{w,p_n}(X_{i,t}) - q_t^{w,\text{new}}(X_{i,t})|^2 \right)$$

$$+ \left( \frac{72}{n_t} \sum_{i=n_l+1}^{n_l+n_t} |\hat{q}_{n,t}^{w,p_n}(X_{i,t}) - q_t^{w,\text{new}}(X_{i,t})|^2 \right.$$

$$- \frac{144}{n_t} \sum_{i=n_l+1}^{n_l+n_t} |\hat{q}_{n,t}^{w,p_n}(X_{i,t}) - q_t(X_{i,t})|^2 \Bigg)$$

$$+ \frac{144}{n_t} \sum_{i=n_l+1}^{n_l+n_t} |\hat{q}_{n,t}^{w,p_n}(X_{i,t}) - q_t(X_{i,t})|^2$$

$$= \sum_{j=1}^{5} T_{j,n}.$$

The proof is completed once we have shown that

$$(6.10) \qquad \mathbf{E} T_{j,n} \leq const \cdot C^{2d/(2p+d)} \cdot \left( \frac{\log n}{n} \right)^{2p/(2p+d)}$$

for $j \in \{1, 2, \ldots, 5\}$.

To apply Bernstein's inequality, we first bound the variance

$$\sigma^2 = \mathbf{V}(|\hat{q}_{n,t}^{w,p}(X_{n_l+1,t}) - q_t(X_{n_l+1,t})|^2 | X_{i,t:w_{\max}(t)+1}^{t,\text{new}} \ (i = 1, \ldots, n_l), \mathcal{D}_{n,t+1}^T)$$

$$\leq \mathbf{E}(|\hat{q}_{n,t}^{w,p}(X_{n_l+1,t}) - q_t(X_{n_l+1,t})|^4 | X_{i,t:w_{\max}(t)+1}^{t,\text{new}} \ (i = 1, \ldots, n_l), \mathcal{D}_{n,t+1}^T)$$

$$\leq 4L^2 \mathbf{E}(|\hat{q}_{n,t}^{w,p}(X_{n_l+1,t}) - q_t(X_{n_l+1,t})|^2 | X_{i,t:w_{\max}(t)+1}^{t,\text{new}} \ (i = 1, \ldots, n_l), \mathcal{D}_{n,t+1}^T)$$

$$= 4L^2 \int |\hat{q}_{n,t}^{w,p}(x) - q_t(x)|^2 \mu_t(dx).$$

Then, because $\hat{q}_{n,t}^w$ and $q_t$ are bounded in absolute value by $L$, we obtain from Bernstein's inequality (see, e.g., Lemma A.2 in [14])

$$\mathbf{P}\{T_{1,n} > \epsilon | X_{i,t:w_{\max}(t)+1}^{t,\text{new}} \ (i = 1, \ldots, n_l), \mathcal{D}_{n,t+1}^T\}$$

$$\leq |\mathcal{P}_n| \cdot \max_{p \in \mathcal{P}_n} \mathbf{P}\Bigg\{ \int |\hat{q}_{n,t}^{w,p}(x) - q_t(x)|^2 \mu_t(dx)$$

$$- \frac{2}{n_t} \sum_{i=n_l+1}^{n_l+n_t} |\hat{q}_{n,t}^{w,p}(X_{i,t}) - q_t(X_{i,t})|^2 > \epsilon$$

$$\Bigg| X_{i,t:w_{\max}(t)+1}^{t,\text{new}} \ (i = 1, \ldots, n_l), \mathcal{D}_{n,t+1}^T \Bigg\}$$

$$= |\mathcal{P}_n| \cdot \max_{p \in \mathcal{P}_n} \mathbf{P}\Bigg\{ \int |\hat{q}_{n,t}^{w,p}(x) - q_t(x)|^2 \mu_t(dx)$$

$$- \frac{1}{n_t} \sum_{i=n_l+1}^{n_l+n_t} |\hat{q}_{n,t}^{w,p}(X_{i,t}) - q_t(X_{i,t})|^2$$

$$> \frac{\epsilon}{2} + \frac{1}{2} \int |\hat{q}_{n,t}^{w,p}(x) - q_t(x)|^2 \mu_t(dx)$$

$$\left| X_{i,t:w_{\max}(t)+1}^{t,\text{new}}(i=1,\ldots,n_l), \mathcal{D}_{n,t+1}^T \right\}$$

$$\leq |\mathcal{P}_n| \cdot \max_{p \in \mathcal{P}_n} \mathbf{P} \left\{ \int |\hat{q}_{n,t}^{w,p}(x) - q_t(x)|^2 \mu_t(dx) \right.$$

$$- \frac{1}{n_t} \sum_{i=n_l+1}^{n_l+n_t} |\hat{q}_{n,t}^{w,p}(X_{i,t}) - q_t(X_{i,t})|^2$$

$$\left. > \frac{\epsilon}{2} + \frac{1}{2} \cdot \frac{\sigma^2}{4L^2} \right| X_{i,t:w_{\max}(t)+1}^{t,\text{new}} \ (i=1,\ldots,n_l), \mathcal{D}_{n,t+1}^T \right\}$$

$$\leq |\mathcal{P}_n| \cdot \exp\left( - \frac{n_t(\epsilon/2 + \sigma^2/(8L^2))^2}{2\sigma^2 + 2(\epsilon/2 + \sigma^2/(8L^2)) \cdot (4L^2/3)} \right)$$

$$\leq |\mathcal{P}_n| \cdot \exp\left( - \frac{n_t(\epsilon/2 + \sigma^2/(8L^2))^2}{(16L^2 + 8L^2/3)(\epsilon/2 + \sigma^2/(8L^2))} \right)$$

$$\leq |\mathcal{P}_n| \cdot \exp\left( - \frac{1}{32 + 16/3} \cdot \frac{n_t \epsilon}{L^2} \right) = |\mathcal{P}_n| \cdot \exp\left( - \frac{3}{112} \cdot \frac{n_t \epsilon}{L^2} \right).$$

Thus,

$$\mathbf{E} T_{1,n} \leq \int_0^\infty \mathbf{P}\{T_{1,n} > s\}\, ds$$

$$= \int_0^\infty \mathbf{E}\{\mathbf{P}\{T_{1,n} > s | X_{i,t:w_{\max}(t)+1}^{t,\text{new}} \ (i=1,\ldots,n_l), \mathcal{D}_{n,t+1}^T\}\}\, ds$$

$$\leq |\mathcal{P}_n| \cdot \int_0^\infty \exp\left( - \frac{3n_t}{112L^2} \cdot s \right) ds$$

$$\leq \frac{112L^2}{3} \cdot \frac{|\mathcal{P}_n|}{n_t} \leq const \cdot C^{2d/(2p+d)} \cdot \left( \frac{\log n}{n} \right)^{2p/(2p+d)}.$$

Furthermore, by $a^2 = (a - b + b)^2 \leq 2(a-b)^2 + 2b^2$, we get

$$T_{2,n} \leq \frac{4}{n_t} \sum_{i=n_l+1}^{n_l+n_t} |q_t^{w,\text{new}}(X_{i,t}) - q_t(X_{i,t})|^2,$$

from which we conclude, together with (3.7) and (6.8), that

$$\mathbf{E} T_{2,n} = \mathbf{E}\{\mathbf{E}\{T_{2,n} | X_{i,t:w_{\max}(t)+1}^{t,\text{new}} \ (i=1,\ldots,n_l), \mathcal{D}_{n,t+1}^T\}\}$$

$$\leq 4\mathbf{E} \int |q_t^{w,\text{new}}(x) - q_t(x)|^2 \mu_t(dx)$$

$$\leq const \cdot C^{2d/(2p+d)} \cdot \left(\frac{\log n}{n}\right)^{2p/(2p+d)}.$$

Similarly, we get

$$\mathbf{E}T_{4,n} \leq const \cdot C^{2d/(2p+d)} \cdot \left(\frac{\log n}{n}\right)^{2p/(2p+d)}.$$

To bound $T_{3,n}$, we apply Lemma 6.1, which shows

$$\mathbf{P}\{T_{3,n} > \epsilon | X_{i,t:w_{\max}(t)+1}^{t,\text{new}} \ (i=1,\dots,n_l), \mathcal{D}_{n,t+1}^T\}$$

$$\leq \mathbf{P}\Bigg\{\frac{1}{n_t}\sum_{i=n_l+1}^{n_l+n_t} |\hat{q}_{n,t}^w(X_{i,t}) - q_t^{w,\text{new}}(X_{i,t})|^2$$

$$> \frac{\epsilon}{4} + 18 \cdot \min_{p\in\mathcal{P}_n} \frac{1}{n_t}\sum_{i=n_l+1}^{n_l+n_t} |\hat{q}_{n,t}^p(X_{i,t}) - q_t^{w,\text{new}}(X_{i,t})|^2$$

$$\Bigg| X_{i,t:w_{\max}(t)+1}^{t,\text{new}} \ (i=1,\dots,n_l), \mathcal{D}_{n,t+1}^T \Bigg\}$$

$$\leq c_5 \cdot \frac{|\mathcal{P}_n|}{\epsilon \cdot n_t}.$$

This implies for any $u > 0$

$$\mathbf{E}T_{3,n} \leq \int_0^\infty \mathbf{P}\{T_{3,n} > \epsilon\}\, d\epsilon$$

$$\leq \int_0^\infty \mathbf{E}\{\mathbf{P}\{T_{3,n} > \epsilon | X_{i,t:w_{\max}(t)+1}^{t,\text{new}} \ (i=1,\dots,n_l), \mathcal{D}_{n,t+1}^T\}\}\, d\epsilon$$

$$\leq u + \int_u^{const} c_5 \cdot \frac{|\mathcal{P}_n|}{\epsilon \cdot n_t}\, d\epsilon$$

$$= u + c_5 \cdot \frac{|\mathcal{P}_n|}{n_t} \cdot (\log(const) - \log u),$$

where we have used that (3.16) and the boundedness of $q_t^{w,\text{new}}$ (which is a consequence of the boundedness of $f_t$ on $[-A, A]^d$) yield

$$T_{3,n} \leq \frac{4}{n_t}\sum_{i=n_l+1}^{n_l+n_t} |\hat{q}_{n,t}^w(X_{i,t}) - q_t^{w,\text{new}}(X_{i,t})|^2 \leq const.$$

With $u = \log(n)/n$, we get

$$\mathbf{E}T_{3,n} \leq \frac{\log n}{n}\left(1 + c_6\left(\log(const) - \log\left(\frac{\log n}{n}\right)\right)\right)$$

$$\leq const \cdot C^{2d/(2p+d)} \cdot \left(\frac{\log n}{n}\right)^{2p/(2p+d)}.$$

Furthermore,

$$(6.11) \quad \begin{aligned} \mathbf{E}T_{5,n} &= \mathbf{E}\{\mathbf{E}\{T_{5,n}|X_{i,t:w_{\max}(t)+1}^{t,\mathrm{new}} \ (i=1,\dots,n_l), \mathcal{D}_{n,t+1}^T\}\} \\ &= 144 \cdot \mathbf{E} \int |\hat{q}_{n,t}^{w,p_n}(x) - q_t(x)|^2 \mu_t(dx). \end{aligned}$$

Consequently, it remains to verify that

$$(6.12) \quad \mathbf{E} \int |\hat{q}_{n,t}^{w,p_n}(x) - q_t(x)|^2 \mu_t(dx) \le const \cdot C^{2d/(2p+d)} \cdot \left(\frac{\log n}{n}\right)^{2p/(2p+d)}$$

for some suitably selected $p_n \in \mathcal{P}_n$.

To bound $\mathbf{E}T_{5,n}$, we use the error decomposition

$$\begin{aligned} \int |\hat{q}&_{n,t}^{w,p_n}(x) - q_t(x)|^2 \mu_t(dx) \\ &= \int |\hat{q}_{n,t}^{w,p_n}(x) - q_t(x)|^2 \mu_t(dx) - \frac{2}{n_l}\sum_{i=1}^{n_l} |\hat{q}_{n,t}^{w,p_n}(X_{i,t}) - q_t(X_{i,t})|^2 \\ &\quad + \frac{2}{n_l}\sum_{i=1}^{n_l} |\hat{q}_{n,t}^{w,p_n}(X_{i,t}) - q_t(X_{i,t})|^2 - \frac{2}{n_l}\sum_{i=1}^{n_l} |\tilde{q}_{n,t}^{w,p_n}(X_{i,t}) - q_t(X_{i,t})|^2 \\ &\quad + \frac{2}{n_l}\sum_{i=1}^{n_l} |\tilde{q}_{n,t}^{w,p_n}(X_{i,t}) - q_t(X_{i,t})|^2 - \frac{4}{n_l}\sum_{i=1}^{n_l} |\tilde{q}_{n,t}^{w,p_n}(X_{i,t}) - q_t^{w,\mathrm{new}}(X_{i,t})|^2 \\ &\quad + \frac{4}{n_l}\sum_{i=1}^{n_l} |\tilde{q}_{n,t}^{w,p_n}(X_{i,t}) - q_t^{w,\mathrm{new}}(X_{i,t})|^2 \\ &= \sum_{j=6}^{9} T_{j,n}, \end{aligned}$$

with

$$p_n = (k, 2^l) \qquad \text{where } l = \lceil \log_2(C^{-2/(2p+d)}(n/\log(n))^{-1/(2p+d)}) \rceil.$$

Because $q_t$ is bounded in absolute value by $L$, we have

$$T_{7,n} \le 0 \quad \text{and} \quad \mathbf{E}T_{7,n} \le 0.$$

Furthermore, in the same way as for $T_{2,n}$, we obtain from (3.7) and (6.8)

$$\begin{aligned} \mathbf{E}T_{8,n} &\le 4\mathbf{E}\left\{\mathbf{E}\left\{\frac{1}{n_l}\sum_{i=1}^{n_l} |q_t(X_{i,t}) - q_t^{w,\mathrm{new}}(X_{i,t})|^2 \Big| \mathcal{D}_{n,t+1}^T\right\}\right\} \\ &= 4\mathbf{E} \int |q_t(x) - q_t^{w,\mathrm{new}}(x)|^2 \mu_t(dx) \\ &\le const \cdot C^{2d/(2p+d)} \cdot \left(\frac{\log n}{n}\right)^{2p/(2p+d)}, \end{aligned}$$

where the last equality follows from the fact that the conditional expectation $q_t^{w,\text{new}}(x)$ does not depend on data from time $t$.

Next, we bound $T_{6,n}$. The functions $\hat{q}_{n,t}^{w,p_n}$ and $q_t$ are bounded in absolute value by $L$, and $\tilde{q}_{n,t}^{w,p_n}$ belongs to the linear vector space $\mathcal{H}_{n,p_n}$, whose dimension $D_n$ is bounded by some constant (depending on $A$ and $k$) times $C^{2d/(2p+d)} \cdot (n/\log(n))^{d/(2p+d)}$. As in the proof of Theorem 11.3 in [14] [in particular, the proof of inequality (11.6)], this implies

$$\mathbf{E}T_{6,n} = \mathbf{E}\{\mathbf{E}\{T_{6,n}|\mathcal{D}_{n,t+1}^T\}\}$$

$$\leq c_7 L^2 \frac{(\log n_l + 1) \cdot C^{2d/(2p+d)} \cdot (n/\log(n))^{d/(2p+d)}}{n_l}$$

$$\leq const \cdot C^{2d/(2p+d)} \cdot \left(\frac{\log n}{n}\right)^{2p/(2p+d)}.$$

Finally, we bound $T_{9,n}$. With

$$\sigma^2 = \sup_{x \in \mathbb{R}^d} \mathbf{E}^*\{|\hat{Y}_{1,t}^{w,\text{new}}|^2|X_{1,t} = x\} \leq 4L^2 < \infty,$$

we can conclude from Theorem 11.1 in [14] that

$$\mathbf{E}\{T_{9,n}|X_{i,t} \ (i=1,\ldots,n_l), \mathcal{D}_{n,t+1}^T\}$$

$$\leq 4\sigma^2 \cdot \frac{D_n}{n_l} + 4 \min_{h \in \mathcal{H}_{n,p_n}} \frac{1}{n_l} \sum_{i=1}^{n_l} |h(X_{i,t}) - q_t^{w,\text{new}}(X_{i,t})|^2$$

$$\leq 4\sigma^2 \cdot C^{2d/(2p+d)} \cdot \frac{c_8}{n^{2p/(2p+d)} \cdot \log(n)^{d/(2p+d)}}$$

$$+ 4 \min_{h \in \mathcal{H}_{n,p_n}} \frac{1}{n_l} \sum_{i=1}^{n_l} |h(X_{i,t}) - q_t^{w,\text{new}}(X_{i,t})|^2.$$

Therefore,

$$\mathbf{E}T_{9,n} = \mathbf{E}\{\mathbf{E}\{T_{9,n}|X_{i,t} \ (i=1,\ldots,n_l), \mathcal{D}_{n,t+1}^T\}\}$$

$$\leq 12\sigma^2 \cdot C^{2d/(2p+d)} \cdot \left(\frac{\log n}{n}\right)^{2p/(2p+d)}$$

$$+ 4 \min_{h \in \mathcal{H}_{n,p_n}} \mathbf{E} \int |h(x) - q_t^{w,\text{new}}(x)|^2 \mu_t(dx)$$

$$\leq 12\sigma^2 \cdot C^{2d/(2p+d)} \cdot \left(\frac{\log n}{n}\right)^{2p/(2p+d)}$$

$$+ 8\mathbf{E} \int |q_t^{w,\text{new}}(x) - q_t(x)|^2 \mu_t(dx)$$

$$+ 8 \min_{h \in \mathcal{H}_{n,p_n}} \int |h(x) - q_t(x)|^2 \mu_t(dx).$$

Note that for the last term in the last inequality (without the factor 8) we get

$$\min_{h \in \mathcal{H}_{n,p}} \int |h(x) - q_t(x)|^2 \mu_t(dx) \leq \min_{h \in \mathcal{H}_{n,p}} \sup_{x \in [-A,A]^d} |h(x) - q_t(x)|^2.$$

Because we have assumed that $q_t$ is $(p, C)$-smooth, there exist a $h \in \mathcal{H}_{n,p}$ with

$$\sup_{x \in [-A,A]^d} |h(x) - q_t(x)| \leq c_9 \cdot C \cdot \delta_n^p,$$

where $\delta_n = C^{-2/(2p+d)} \cdot (n/\log(n))^{-1/(2p+d)}$ is the edge length in the cubic partition used in the definition of the spline space; see Theorem 12.8 in [24]. We conclude that

$$\min_{h \in \mathcal{H}_{n,p}} \int |h(x) - q_t(x)|^2 \mu_t(dx)$$
$$\leq c_9^2 \cdot C^2 \cdot \delta_n^{2p}$$
$$= c_9^2 \cdot C^2 \cdot C^{-4p/(2p+d)} \cdot (n/\log(n))^{-2p/(2p+d)}$$
$$\leq const \cdot C^{2d/(2p+d)} \cdot \left(\frac{\log n}{n}\right)^{2p/(2p+d)}.$$

From (3.7), (6.8) and the above inequality we see that

$$\mathbf{E} T_{9,n} \leq const \cdot C^{2d/(2p+d)} \cdot \left(\frac{\log n}{n}\right)^{2p/(2p+d)}$$

has an upper bound with the proper rate. The proof of Theorem 4.4 is complete.  □

**Acknowledgments.** The authors would like to thank the associate editor and the referees for the very constructive suggestions which helped us to improve the article.

D. Egloff
Zurich Cantonal Bank
P.O. Box, CH-8010 Zurich
Switzerland
E-mail: daniel.egloff@zkb.ch

M. Kohler
N. Todorovic
Department of Mathematics
University of Saarbrücken
Postfach 151150, D-66041 Saarbrücken
Germany
E-mail: kohler@math.uni-sb.de
        todorovic@math.uni-sb.de